\documentclass[11pt]{amsart}

\usepackage{amsmath,amsthm, amscd, amssymb, amsfonts}
\usepackage[all]{xy}

\usepackage[dvips, dvipsnames, usenames]{color}

\newtheorem{theorem}{Theorem}[section]

\newtheorem{coro}[theorem]{Corollary}

\newtheorem{prop}[theorem]{Proposition}

\theoremstyle{definition}
\newtheorem{definition}[theorem]{Definition}
\newtheorem{example}[theorem]{Example}

\theoremstyle{remark}
\newtheorem{remark}[theorem]{Remark}

\newcommand{\ydh}{{}^{H}_{H}\mathcal{YD}}

\def\zt{\Z^{\theta}}

\def\G{\mathbb{G}}
\newcommand\id{\operatorname{id}}

\newcommand\ord{\operatorname{ord}}
\newcommand\Hom{\operatorname{Hom}}

\newcommand\Aut{\operatorname{Aut}}

\newcommand\ad{\operatorname{ad}}
\newcommand\sop{\operatorname{sop}}

\newcommand{\ku}{ \mathbf{k}}

\def\ot{\otimes}

\def\Z{\mathbb{Z}}
\def\N{\mathbb{N}}
\def\xx{\mathbb{X}}

\def\cB{\mathcal{B}}
\def\cO{\mathcal{O}}
\def\cP{\mathcal{P}}
\def\cC{\mathcal{C}}
\def\cR{\mathcal{R}}
\def\cW{\mathcal{W}}
\def\Xf{\mathcal{X}}
\def\cU{\mathcal{U}}
\def\qmb{\mathbf{q}}

\newcommand{\vi}{\textbf{(i)} }
\newcommand{\vii}{\textbf{(ii)} }
\newcommand{\viii}{\textbf{(iii)} }
\newcommand{\viv}{\textbf{(iv)} }
\newcommand{\vv}{\textbf{(v)} }
\newcommand{\vvi}{\textbf{(vi)} }
\newcommand{\vvii}{\textbf{(vii)} }
\newcommand{\vviii}{\textbf{(viii)} }
\newcommand{\vix}{\textbf{(ix)} }

\newcommand{\schi}{{s_p^*\chi}}

\def\pf{\begin{proof}}
\def\epf{\end{proof}}

\begin{document}


\title[Nichols algebras of unidentified diagonal type]{Nichols algebras of unidentified diagonal type}
\author[Iv\'an Angiono]{Iv\'an Angiono}

\address{FaMAF-CIEM (CONICET), Universidad Nacional de C\'ordoba,
Medina A\-llen\-de s/n, Ciudad Universitaria, 5000 C\' ordoba, Rep\'
ublica Argentina.} \email{angiono@famaf.unc.edu.ar}

\thanks{The work was partially supported by CONICET,
FONCyT-ANPCyT, Secyt (UNC), Mincyt (C\'ordoba)}

\begin{abstract}
The Nichols algebras of diagonal type with finite root system are
either of standard, super or (yet) unidentified type. A concrete
description of the defining relations of all those Nichols algebras
was given in \cite{A-exp presentation}. In the present paper we use
this result to give an explicit presentation of all Nichols algebras
of unidentified type.
\end{abstract}

\maketitle

\section*{Introduction}

The classification of braidings of diagonal type whose Nichols
algebra has a finite number of roots was given in \cite{H-classif
RS}. This problem is related with the classification of
finite-dimensional pointed Hopf algebras over abelian groups. The
list of Heckenberger can be split off in three
families:
\begin{itemize}
  \item standard braidings, introduced in \cite{AA};
  \item braidings of super type, see \cite{AAY};
  \item a finite list of braidings whose connected components have rank
  less than eight. We call them \emph{unidentified}.
\end{itemize}
We recall that the Nichols algebra of a braided vector space $(V,c)$
is a quotient of its tensor algebra by a suitable ideal $I(V)$. A
crucial question involving these Nichols algebras is to obtain a
minimal set of relations generating $I(V)$. Such problem was solved
in \cite{A-Standard} for the first family, with a formula for the
dimension of each Nichols algebra of standard type. A presentation
for the second family was given in \cite{Y} for the generic case,
and in \cite{AAY} for the non-generic case (except by some
considerations for small orders on the entries of the braiding
matrix).

A complete answer can be found in \cite{A-exp presentation}, where
the main result gives a list of relations satisfied by the
generators of the Nichols algebras, depending on the matrix entries,
see Theorem \ref{thm:presentacion minima} below. This paper depends
strongly on \cite{A-presentation}, the key point to obtain the
desired presentation. It remains the problem to identify the
relations needed for each one of the braidings. This is instrumental
for several important questions concerning pointed Hopf algebras,
among them the explicit determination of all liftings and their
representation theory.

In this paper, we deal with unidentified braidings and give a complete list of relations
generating the defining ideal for the Nichols algebra of each braiding of this kind.
For the small ranks we give also the list of positive roots for each
case and the dimension.

The organization of the paper is the following. The first Section
includes some notions about PBW bases and root systems of Nichols
algebras of diagonal type, as well as the presentation by generators
and relations from \cite{A-exp presentation}. The second Section is
devoted to give explicitly the presentation for each diagram. We
begin with some unidentified families in rank 2, 3 or 4, and finally
we split the remaining braidings in families by considering certain
similarities on the associated generalized Cartan matrices.

\subsection*{Notation} $\ku$ will denote an algebraically closed field of characteristic zero. For each $N>1$, $\G_N$ will denote
the group of $N$-roots of unity in $\ku$, and $\G'_N$ the corresponding subset of primitive roots of order $N$.

\section{Preliminaries}

We will recall all the preliminary results and fix the notation that
we will use along this work. They are related with the theory of PBW
bases for braided Hopf algebras of diagonal type \cite{Kh}, the Weyl
groupoid of diagonal braided vector spaces \cite{H-Weyl gpd} and the
presentation of Nichols algebras by generators and relations
\cite{A-exp presentation}.

\subsection{PBW bases for Nichols algebras of diagonal type}\label{subsection:pbw}

\

We begin with the definition of a Nichols algebra associated to a braided vector space. To this end, fix
a Hopf algebra $H$ with bijective antipode, and denote by $\ydh$ the category of left Yetter-Drinfeld modules over
$H$. Let $V\in\ydh$. The tensor algebra $T(V)$ admits a braiding extending $c:V\ot V\to V\ot V$, and under such braiding it has a unique
structure of graded braided Hopf algebra in $\ydh$ such that $V \subseteq \cP(V)$ (i.e. $\Delta(x)=x\ot 1+1\ot x$).

\begin{definition}{\cite{AS Pointed HA}} \label{def:algebra de Nichols}
Let $\mathfrak{S}$ be the family of all the homogeneous Hopf ideals of $I\subseteq T(V)$ such that
\begin{itemize}
    \item $I$ is generated by homogeneous elements of degree $\geq 2$,
    \item $I$ is a Yetter-Drinfeld submodule of $T(V)$.
\end{itemize}
The \emph{Nichols algebra} $\cB(V)$ associated to $V$ is the quotient of $T(V)$ by the biggest ideal $I(V)$ of $\mathfrak{S}$.
\end{definition}
\medskip

The definition does not depend on the realization of a braided vector space $(V,c)$ as a Yetter-Drinfeld module.
In particular we will consider braidings of diagonal type; that is, there exists a basis $\{x_i\}_{i\in I}$ of $V$
and a family of non-zero scalars $(q_{ij})_{i,j\in I}$ such that $c(x_i\ot x_j)=q_{ij}\, x_j\ot x_i$. These braided
vector spaces are related with vector spaces over group algebras of finite abelian groups.

Fix $(V,c)$ a braided vector space of diagonal type. We will describe, following \cite{Kh}, a particular PBW basis for
each graded braided Hopf algebra $\cB=\oplus_{n\in \N} \cB^n$ generated by $\cB^1 \cong V$ as an algebra.

Assume that $V$ is finite dimensional, and denote $\theta=\dim V$. Fix a basis $X= \{x_1,\dots, x_{\theta}\}$ of $V$ as above.
$\xx$ will denote the set of words with letters in $X$. We consider the lexicographical order on $\xx$. We can identify
$\ku \xx$ with $T(V)$.

$T(V)$ admits a unique $\Z^I$-graduation as a braided Hopf algebra such that $\deg x_i=\alpha_i$, where $(\alpha_i)_{i\in I}$
is the canonical basis of $\Z^I$. Assume that $\dim V=\theta<\infty$.
Let $\chi: \zt\times \zt \to \ku^{\times}$ be the bicharacter determined by the condition
\begin{equation}\label{eq:forma bilineal diagonal}
\chi(\alpha_i, \alpha_j) = q_{ij}, \quad \mbox{ for each pair }1\le i, j \le \theta.
\end{equation}
Then, for each pair of $\zt$-homogeneous elements $u,v \in \xx$,
\begin{equation}\label{eq:braiding tipo diagonal}
    c(u \ot v)= q_{u,v} v \ot u, \qquad q_{u,v} = \chi(\deg u, \deg v)\in \ku^{\times}.
\end{equation}

Note that \cite[Prop. 2.10]{AS Pointed HA} implies that $I(V)$ is a $\zt$-homogeneous ideal,
and then $\cB(V)$ is $\zt$-graded, see also \cite[Prop. 1.2.3]{L}.

\bigskip

\begin{definition}
$u \in \xx -\{ 1\}$ is a \emph{Lyndon word} if for every decomposition $u=vw$, $v,w \in\xx - \left\{ 1 \right\}$, it holds that $u<w$.
We will denote the set of all Lyndon words by $L$.
\end{definition}

We know that any word $u \in \xx$ admits a unique decomposition as a non-increasing product of Lyndon words:
\begin{equation}\label{eq:descly}
u=l_1l_2\dots  l_r, \qquad l_i \in L, l_r \leq \dots \leq l_1.
\end{equation}
It is called the \emph{Lyndon decomposition} of $u \in \xx$, and each $l_i\in L$ in \eqref{eq:descly} is called a
\emph{Lyndon letter} of $u$.

For each $u \in L-X$, the \emph{Shirshov decomposition} of $u$ is the decomposition $u=u_1u_2$, $u_1,u_2\in L$, such
that $u_2$ is the smallest end of $u$ between all the possible decompositions satisfying these conditions (it is
easily proved that each Lyndon word admits at least one of such decompositions).
\medskip

\medskip

For a general braided vector space, the \emph{braided bracket} of $x,y\in T(V)$ is defined by
\begin{equation}\label{eq:braidedcommutator}
[x,y]_c := \text{multiplication } \circ \left( \id - c \right) \left( x \ot y \right).
\end{equation}

Using the previous decompositions and the identification of $\xx$ with a basis of $T(V)$, we can define a $\ku$-linear
endormorphism $\left[ - \right]_c$ of $T(V)$ as follows:
$$ \left[ u \right]_c := \begin{cases} u,& \text{if } u = 1 \text{ or }u \in X;\\
[\left[ v \right]_c, \left[ w \right]_c]_c,  & \text{if } u \in L, \, \ell(u)>1, \ u=vw \text{ is the Shirshov decomposition};\\
\left[u_1\right]_c \dots \left[u_t\right]_c,& \text{ if } u\in \xx-L \text{ and its Lyndon decomposition is }u=u_1\dots u_t.
\end{cases}
$$
We will obtain PBW bases using this automorphism.

\begin{definition} The \emph{hyperletter} corresponding to $l \in L$ is the element $\left[l\right]_c$. An \emph{hyperword} is a
word written in hyperletters; a \emph{monotone hyperword} is an hyperword $W=\left[u_1\right]_c^{k_1}\dots\left[u_m\right]_c^{k_m}$
such that $u_1>\dots >u_m$.
\end{definition}

A different order on $\xx$, considered in \cite{U} and used implicitly in \cite{Kh} is the the \emph{deg-lex order}, defined as
follows. For each pair $u,v \in \xx$, we say that $u \succ v$ if $\ell(u)<\ell(v)$, or
$\ell(u)=\ell(v)$ and $u>v$ for the lexicographical order. Such order is total. The empty word $1$ is the maximal element for
$\succ$, and this order is invariant by left and right multiplication.
\medskip

\noindent In what follows $I$ will denote a Hopf ideal, and $R=T(V)/I$. Let $\pi: T(V) \rightarrow R$ be the canonical
projection. We set:
$$G_I:= \left\{ u \in \xx: u \notin \\ \ku \xx_{\succ u}+I  \right\}.$$
Therefore, if $u \in G_I$ and $u=vw$, then $v,w \in G_I$. In this way, each $u \in G_I$ is a non-increasing product of
Lyndon words of $G_I$.

Consider the set $S_I:=G_I\cap L$, and define $h_I:S_I\to\left\{2,3,\dots \right\}\cup \{\infty\}$ by the condition:
\begin{equation}\label{defaltura}
    h_I(u):= \min \left\{ t \in \N : u^t  \in \ku \xx_{\succ u^t} + I \right\}.
\end{equation}
Following \cite{Kh} we have the following results.

\begin{theorem}\label{thm:base PBW Kharchenko} The set
$$ \{ \left[u_1\right]_c^{k_1}\dots\left[u_m\right]_c^{k_m}: \, m\in\N_0, u_1> \ldots >u_m, u_i\in S_I, 0<k_i<h_I(u_i) \} $$
is a PBW basis of $H=T(V)/I$.
\qed
\end{theorem}

\begin{coro}\label{cor:primero}
\vi A word $u$ belongs to $G_I$ if and only if the hyperletter $\left[u\right]_c$ is not a linear combination of greater hyperwords
$\left[ w \right]_c$, $w \succ u$, whose hyperletters are in $S_I$, modulo $I$. \qed
\medskip

\noindent \vii If $v\in S_I$ is such that $h_I(v)<\infty$, then $q_{v,v}$ is a root of unity. Moreover, if $\ord q_{v,v}=h$, then $h_I(v)=h$, and $\left[v\right]^h$ is a
linear combination of hyperwords $\left[w\right]_c$, $w\succ v^h$.
\qed
\end{coro}
\medskip

Let $\Delta^V_+$ be the set of degrees of a PBW basis of $\cB(V)$, counted with their multiplicities \cite{H-Weyl gpd}. We can see
that it does not depend on the PBW basis, \cite{H-Weyl gpd,AA}. We can attach a Cartan scheme $\cC$, a Weyl groupoid $\cW$ and a
root system $\cR$ in the sense of \cite{CH1,HY}, see \cite[Thms. 6.2, 6.9]{HS}. To this end, define for each
$1 \leq i \neq j \leq \theta$,
\begin{equation}
-a_{ij}:= \min \left\{ n \in \mathbb{N}_0: (n+1)_{q_{ii}} (1-q_{ii}^n\widetilde{q_{ij}} )=0 \right\}, \label{defn:mij}
\end{equation}
and set $a_{ii}=2$. The symmetry $s_i\in\Aut(\Z^\theta)$ is defined by the condition $s_i(\alpha_j)=\alpha_j-a_{ij}\alpha_i$.

Set $\widetilde{q}_{rs}=\chi(s_i(\alpha_r),s_i(\alpha_s))$. Let $V_i$ be another vector space of the same dimension, and attach
to it the matrix
$\widetilde{\qmb}=(\widetilde{q}_{rs})$. By \cite{H-Weyl gpd},
$$ \Delta^{V_i}_+ = s_i\left( \Delta^V_+ \setminus\{\alpha_i\}\right) \cup \{\alpha_i\}. $$
Therefore last equation lets us to define the Weyl groupoid of $V$, whose root system is defined by the sets
$\Delta^{V'}$, $V'$ obtained after to apply some reflections to the matrix of $V$. Those braided vector spaces obtained after to
apply the symmetries $s_i$ define the \emph{Weyl equivalence class} of $V$.

When the root system is finite, we can prove that each root is real, and in consequence it has multiplicity one, see \cite{CH1}.

\subsection{A presentation by generators and relations of Nichols algebras of diagonal type}\label{subsection:presentation}

Fix as above a finite-dimensional braided vector space $(V,c)$ of diagonal
type, with braiding matrix $(q_{ij})_{1\leq i,j\leq\theta}$,
$\theta=\dim V$, and a basis $x_1,\ldots,x_\theta$ of $V$ such
that $c(x_i\ot x_j)=q_{ij} x_j\ot x_i$. Let $\chi$ be the
bicharacter associated to $(q_{ij})$.

We denote $\widetilde{q_{ij}}=q_{ij}q_{ji}$. The \emph{generalized
Dynkin diagram} of a matrix $(q_{ij})_{1\leq i,j\leq\theta}$ is a
graph with $\theta$ vertices, labeled with the scalars $q_{ii}$, and
an arrow between the vertices $i$ and $j$ if $\widetilde{q_{ij}}\neq
1$, labeled with this scalar. For example, given $q\in\ku^\times$,
the matrices
$$ \left( \begin{array}{ccc} q^2 & q^{-2} & 1 \\ 1 & q^2 & q^{-2} \\ 1 & 1 & q \end{array} \right),
\left( \begin{array}{ccc} q^2 & q^{-1} & q \\ q^{-1} & q^2 & q^{-1}
\\ q^{-1} & q^{-1} & q \end{array} \right)$$
have the diagram: $\xymatrix{\circ^{q^2}\ar@{-}[r]^{q^{-2}} &
\circ^{q^2} \ar@{-}[r]^{q^{-2}} & \circ^{q}}$. In fact, two braided vector spaces of diagonal type are \emph{twist equivalent} \cite{AS Pointed HA} if they have the same
generalized Dynkin diagram.

We denote also
$$ x_{i_1i_2\cdots i_k}=(\ad_c x_{i_1})\cdots(\ad_c x_{i_{k-1}})x_{i_k}, \qquad i_j\in\{1,\ldots,\theta\}. $$
For each $m\in \N$, we define the elements
$x_{(m+1)\alpha_i+m\alpha_j}\in\cU(\chi)$ recursively:
\begin{itemize}
  \item if $m=1$, $x_{2\alpha_i+\alpha_j}:= (\ad_c x_i)^2 x_j=x_{iij}$,
  \item $x_{(m+2)\alpha_i+(m+1)\alpha_j}:= [ x_{(m+1)\alpha_i+m\alpha_j}, x_{ij} ]_c$.
\end{itemize}

Call $x_\alpha$, $\alpha\in\Delta_+^V$, the generator of the Kharchenko's PBW basis. We denote
$$ N_\alpha:=\ord \chi(\alpha,\alpha), \qquad \mbox{if } \chi(\alpha,\alpha) \mbox{ is a root of unity.} $$

We give now the main result of \cite{A-exp presentation}.

\begin{theorem}\label{thm:presentacion minima}
Assume that the root system $\Delta^\chi$ is finite. Then $\cB(V)$ admits a presentation by generators
$x_1,\ldots,x_\theta$ and relations:
\begin{align}
&x_\alpha^{N_\alpha}, &\alpha \in \cO(\chi);\label{eqn:potencia
raices}
\\ &(\ad_cx_i)^{m_{ij}+1}x_j, & q_{ii}^{m_{ij}+1} \neq 1;\label{eqn:relacion quantum Serre}
\\ &x_i^{N_i}, & i \mbox{ is not a Cartan vertex};\label{eqn:potencia raices simples}
\end{align}

\noindent $\odot$ if $i,j \in \{1, \ldots, \theta \}$ are such that
$q_{ii}=\widetilde{q_{ij}}=q_{jj}=-1$, and there exists $k\neq i,j$
such that $\widetilde{q_{ik}}^2\neq1$ or
$\widetilde{q_{jk}}^2\neq1$,
\begin{equation}\label{eqn:relacion dos vertices con -1}
x_{ij}^2;
\end{equation}

\noindent $\odot$ if $i,j,k \in \{1, \ldots, \theta \}$ are such
that $q_{jj}=-1$,
$\widetilde{q_{ik}}=\widetilde{q_{ij}}\widetilde{q_{kj}}=1$,
\begin{equation}\label{eqn:relacion vertice -1}
\left[ x_{ijk} , x_j \right]_c;
\end{equation}

\noindent $\odot$ if $i,j \in \{1, \ldots, \theta \}$ are such that
$q_{jj}=-1$, $q_{ii}\widetilde{q_{ij}}\in \G_6$, and also $q_{ii}\in
\G_3$ or $m_{ij}\geq 3$,
\begin{equation}\label{eqn:relacion estandar B2}
\left[ x_{iij}, x_{ij} \right]_c;
\end{equation}

\noindent $\odot$ if $i,j,k \in \{1, \ldots, \theta \}$ are such
that $q_{ii}=\pm \widetilde{q_{ij}}\in\G_3$, $\widetilde{q_{ik}}=1$,
and also $-q_{jj}=\widetilde{q_{ij}}\widetilde{q_{jk}}=1$ or
$q_{jj}^{-1}=\widetilde{q_{ij}}=\widetilde{q_{jk}}\neq -1$,
\begin{equation}\label{eqn:relacion estandar B3}
\left[ x_{iijk} , x_{ij} \right]_c;
\end{equation}

\noindent $\odot$ if $i,j,k \in \{1, \ldots, \theta \}$ are such
that $\widetilde{q_{ik}}, \widetilde{q_{ij}}, \widetilde{q_{jk}}
\neq 1$,
\begin{equation}\label{eqn:relacion triangulo}
x_{ijk}-\frac{1-\widetilde{q_{jk}}}{q_{kj}(1-\widetilde{q_{ik}})}\left[x_{ik},x_j\right]_c-q_{ij}(1-\widetilde{q_{jk}})
\ x_jx_{ik};
\end{equation}

\noindent $\odot$ if $i,j,k\in\{1,\ldots,\theta\}$ are such that one
of the following situations
\begin{itemize}
 \item[$\circ$] $q_{ii}=q_{jj}=-1$, $\widetilde{q_{ij}}^2= \widetilde{q_{jk}}^{-1}$, $\widetilde{q_{ik}}=1$, or
 \item[$\circ$] $\widetilde{q_{ij}}=q_{jj}=-1$, $q_{ii}= -\widetilde{q_{jk}}^2\in\G_3$, $\widetilde{q_{ik}}=1$, or
 \item[$\circ$] $q_{kk}=\widetilde{q_{jk}}=q_{jj}=-1$, $q_{ii}= -\widetilde{q_{ij}}\in\G_3$, $\widetilde{q_{ik}}=1$, or
 \item[$\circ$] $q_{jj}=-1$, $\widetilde{q_{ij}}=q_{ii}^{-2}$, $\widetilde{q_{jk}}=-q_{ii}^{-3}$, $\widetilde{q_{ik}}=1$, or
 \item[$\circ$] $q_{ii}=q_{jj}=q_{kk}=-1$, $\pm\widetilde{q_{ij}}=\widetilde{q_{jk}}\in\G_3$,
 $\widetilde{q_{ik}}=1$,
\end{itemize}
\begin{equation}\label{eqn:relacion super C3}
\left[ \left[x_{ij}, x_{ijk} \right]_c, x_j \right]_c;
\end{equation}

\noindent $\odot$ if $i,j,k\in\{1,\ldots,\theta\}$ are such that
$q_{ii}=q_{jj}=-1$, $\widetilde{q_{ij}}^3=\widetilde{q_{jk}}^{-1}$,
$\widetilde{q_{ik}}=1$,
\begin{equation}\label{eqn:relacion super G3}
\left[ \left[x_{ij}, \left[x_{ij}, x_{ijk} \right]_c \right]_c, x_j
\right]_c;
\end{equation}

\noindent $\odot$ if $i,j,k\in\{1,\ldots,\theta\}$ are such that
$q_{jj}=\widetilde{q_{ij}}^2=\widetilde{q_{jk}}\in \G_3$,
$\widetilde{q_{ik}}=1$,
\begin{equation}\label{eqn:relacion super C3 raiz de orden 3}
\left[ \left[ x_{ijk} , x_j \right]_c x_j \right]_c;
\end{equation}

\noindent $\odot$ if $i,j,k\in\{1,\ldots,\theta\}$ are such that
$q_{kk}=q_{jj}=\widetilde{q_{ij}}^{-1}=\widetilde{q_{jk}}^{-1}\in
\G_9$, $\widetilde{q_{ik}}=1$, $q_{ii}=q_{kk}^6$
\begin{equation}\label{eqn:relacion fila 18 rango 3}
\left[ \left[ x_{iij} , x_{iijk} \right]_c, x_{ij} \right]_c;
\end{equation}

\noindent $\odot$ if $i,j,k\in\{1,\ldots,\theta\}$ are such that
$q_{ii}=\widetilde{q_{ij}}^{-1}\in \G_9$,
$q_{jj}=\widetilde{q_{jk}}^{-1}=q_{ii}^5$, $\widetilde{q_{ik}}=1$,
$q_{kk}=q_{ii}^6$
\begin{equation}\label{eqn:relacion fila 18 rango 3, caso 2}
[\left[x_{ijk}, x_{j} \right]_c, x_k]_c -(1 +
\widetilde{q_{jk}})^{-1}q_{jk} \left[ \left[x_{ijk}, x_{k} \right]_c
, x_{j} \right]_c;
\end{equation}

\noindent $\odot$ if $i,j,k\in\{1,\ldots,\theta\}$ are such that
$q_{jj}=\widetilde{q_{ij}}^3=\widetilde{q_{jk}}\in \G_4$,
$\widetilde{q_{ik}}=1$,
\begin{equation}\label{eqn:relacion super G3 raiz de orden 4}
\left[ \left[ \left[ x_{ijk} , x_j \right]_c, x_j \right]_c, x_j
\right]_c;
\end{equation}

\noindent $\odot$ if $i,j,k\in\{1,\ldots,\theta\}$ are such that
$q_{ii} = \widetilde{q_{ij}} =-1$, $q_{jj}= \widetilde{q_{jk}}^{-1}
\neq-1$, $\widetilde{q_{ik}}=1$,
\begin{equation}\label{eqn:relacion parecida a super C3}
\left[x_{ij}, x_{ijk} \right]_c;
\end{equation}

\noindent $\odot$ if $i,j,k \in\{1,\ldots,\theta\}$ are such that
$q_{ii}= q_{kk} =-1$, $\widetilde{q_{ik}}=1$, $\widetilde{q_{ij}}
\in \G_3$, $q_{jj}= -\widetilde{q_{jk}} = \pm \widetilde{q_{ij}}$,
\begin{equation}\label{eqn:relacion parecida a super C3-bis}
[x_i, x_{jjk}]_c -(1 + q_{jj}^2)q_{kj}^{-1} \left[x_{ijk}, x_{j}
\right]_c - (1 + q_{jj}^2)(1 + q_{jj}) q_{ij} x_j x_{ijk};
\end{equation}

\noindent $\odot$ if $i,j,k,l\in\{1,\ldots,\theta \}$ are such that
$q_{jj}\widetilde{q_{ij}}= q_{jj}\widetilde{q_{jk}}=1$, $q_{kk}=-1$,
$\widetilde{q_{ik}}=\widetilde{q_{il}}=\widetilde{q_{jl}}=1$,
$\widetilde{q_{jk}}^2= \widetilde{q_{lk}}^{-1}= q_{ll}$,
\begin{equation}\label{eqn:relacion super C4}
\left[\left[\left[x_{ijkl},x_k\right]_c, x_j \right]_c, x_k
\right]_c;
\end{equation}

\noindent $\odot$ if $i,j,k,l\in\{1,\ldots,\theta \}$ are such that
$\widetilde{q_{jk}}= \widetilde{q_{ij}}= q_{jj}^{-1}\in
\G_4'\cup\G_6'$, $q_{ii}=q_{kk}=-1$,
$\widetilde{q_{ik}}=\widetilde{q_{il}}=\widetilde{q_{jl}}=1$,
$\widetilde{q_{jk}}^3= \widetilde{q_{lk}}$,
\begin{equation}\label{eqn:relacion super C4 modificada}
\left[\left[x_{ijk},\left[x_{ijkl}, x_k \right]_c\right]_c, x_{jk}
\right]_c;
\end{equation}

\noindent $\odot$ if $i,j,k,l\in\{1,\ldots,\theta \}$ are such that
$q_{ll}=\widetilde{q_{lk}}^{-1}=
q_{kk}=\widetilde{q_{jk}}^{-1}=q^2$, $\widetilde{q_{ij}}=
q_{ii}^{-1}=q^3$ for some $q\in \ku^\times$, $q_{jj}=-1$,
$\widetilde{q_{ik}}=\widetilde{q_{il}}=\widetilde{q_{jl}}=1$,
\begin{equation}\label{eqn:relacion super F4-1}
\left[\left[\left[x_{ijk},x_j\right]_c, \left[x_{ijkl},x_j\right]_c
\right]_c, x_{jk} \right]_c;
\end{equation}

\noindent $\odot$ if $i,j,k,l\in\{1,\ldots,\theta \}$ are such that
one of the following situations hold
\begin{itemize}
  \item[$\circ$] $q_{kk}=-1$, $q_{ii}=\widetilde{q_{ij}}^{-1}= q_{jj}^2$,
$\widetilde{q_{kl}}= q_{ll}^{-1}= q_{jj}^3$, $\widetilde{q_{jk}}=
q_{jj}^{-1}$, $\widetilde{q_{ik}}=\widetilde{q_{il}}=\widetilde{q_{jl}}=1$, or
  \item[$\circ$] $q_{ii}=\widetilde{q_{ij}}^{-1}= -q_{ll}^{-1}=-\widetilde{q_{kl}}$,
$q_{jj}=\widetilde{q_{jk}}=q_{kk}=-1$,
$\widetilde{q_{ik}}=\widetilde{q_{il}}=\widetilde{q_{jl}}=1$,
\end{itemize}
\begin{equation}\label{eqn:relacion super F4-2}
\left[\left[x_{ijkl}, x_j \right]_c, x_k \right]_c-
q_{jk}(\widetilde{q_{ij}}^{-1}-q_{jj})
\left[\left[x_{ijkl},x_k\right]_c, x_j \right]_c;
\end{equation}

\noindent $\odot$ if $i,j,k\in\{1,\ldots,\theta\}$ are such that
$\widetilde{q_{jk}}=1$,
$q_{ii}=\widetilde{q_{ij}}=-\widetilde{q_{ik}}\in \G_3$,
\begin{equation}
\left[x_i, \left[ x_{ij},x_{ik} \right]_c
\right]_c+q_{jk}q_{ik}q_{ji} \left[ x_{iik} ,x_{ij} \right]_c+q_{ij}
\, x_{ij} x_{iik};\label{eqn:relacion mji,mjk=2}
\end{equation}

\noindent $\odot$ if $i,j,k\in\{1,\ldots,\theta\}$ are such that
$q_{jj}=q_{kk}=\widetilde{q_{jk}}=-1$,
$q_{ii}=-\widetilde{q_{ij}}\in\G_3$, $\widetilde{q_{ik}}=1$,
\begin{equation}\label{eqn:relacion especial rank3}
\left[x_{iijk}, x_{ijk} \right]_c;
\end{equation}

\noindent $\odot$ if $i,j\in\{1,\ldots,\theta\}$ are such that
$-q_{ii}, -q_{jj}, q_{ii}\widetilde{q_{ij}},
q_{jj}\widetilde{q_{ij}} \neq 1$,
\begin{equation}\label{eqn:relacion mij,mji mayor que 1}
(1-\widetilde{q_{ij}})q_{jj}q_{ji}\left[x_i, \left[ x_{ij}, x_j
\right]_c \right]_c -
(1+q_{jj})(1-q_{jj}\widetilde{q_{ij}})x_{ij}^2;
\end{equation}

\noindent $\odot$ if $i,j\in\{1,\ldots,\theta\}$ are such that either $m_{ij}\in \{4,5\}$, or else $q_{jj}=-1$, $m_{ij}=3$, $q_{ii} \in \G_4'$,
\begin{equation}\label{eqn:relacion mij mayor que dos, raiz alta}
\left[x_i,x_{3\alpha_i+2\alpha_j}\right]_c-\frac{1-q_{ii}\widetilde{q_{ij}}-q_{ii}^2\widetilde{q_{ij}}^2q_{jj}}{(1-q_{ii}\widetilde{q_{ij}})q_{ji}}
x_{iij}^2;
\end{equation}

\noindent $\odot$ if $i,j\in\{1,\ldots,\theta\}$ are such that
$4\alpha_i+3\alpha_j\notin \Delta_+^\chi$, $q_{jj}=-1$ or
$m_{ji}\geq2$, and also $m_{ij}\geq 3$, or $m_{ij}=2$,
$q_{ii}\in\G_3$,
\begin{equation}\label{eqn:relacion (m+1)alpha i+m alpha j, caso 2}
x_{4\alpha_i+3\alpha_j}=[x_{3\alpha_i+2\alpha_j}, x_{ij} ]_c;
\end{equation}

\noindent $\odot$ if $i,j\in\{1,\ldots,\theta\}$ are such that
$3\alpha_i+2\alpha_j\in\Delta_+^\chi$,
$5\alpha_i+3\alpha_j\notin\Delta_+^\chi$, and
$q_{ii}^3\widetilde{q_{ij}}, q_{ii}^4\widetilde{q_{ij}}\neq 1$,
\begin{equation}\label{eqn:relacion con 2alpha i+alpha j, caso 2}
[x_{iij}, x_{3\alpha_i+2\alpha_j}]_c;
\end{equation}

\noindent $\odot$ if $i,j\in\{1,\ldots,\theta\}$ are such that
$4\alpha_i+3\alpha_j\in\Delta_+^\chi$,
$5\alpha_i+4\alpha_j\notin\Delta_+^\chi$,
\begin{equation}\label{eqn:relacion (m+1)alpha i+m alpha j, caso 3}
x_{5\alpha_i+4\alpha_j}=[x_{4\alpha_i+3\alpha_j}, x_{ij} ]_c;
\end{equation}

\noindent $\odot$ if $i,j\in\{1,\ldots,\theta\}$ are such that
$5\alpha_i+2\alpha_j\in\Delta_+^\chi$, $7\alpha_i+3\alpha_j \notin
\Delta_+^\chi$,
\begin{equation}\label{eqn:relacion con 2alpha i+alpha j, caso 1}
[[x_{iiij}, x_{iij}], x_{iij} ]_c;
\end{equation}

\noindent $\odot$ if $i,j\in\{1,\ldots,\theta\}$ are such that
$q_{jj}=-1$, $5\alpha_i+4\alpha_j\in\Delta_+^\chi$,
\begin{equation}\label{eqn:relacion potencia alta}
[x_{iij},x_{4\alpha_i+3\alpha_j}]_c-
\frac{b-(1+q_{ii})(1-q_{ii}\zeta)(1+\zeta+q_{ii}\zeta^2)q_{ii}^6\zeta^4}
{a\ q_{ii}^3q_{ij}^2q_{ji}^3} x_{3\alpha_i+2\alpha_j}^2,
\end{equation}
where $\zeta=\widetilde{q_{ij}}$,
$a=(1-\zeta)(1-q_{ii}^4\zeta^3)-(1-q_{ii}\zeta)(1+q_{ii})q_{ii}\zeta$,
$b=(1-\zeta)(1-q_{ii}^6\zeta^5)-a\ q_{ii}\zeta$.
\qed
\end{theorem}
\bigskip

\section{Unidentified Nichols algebras of rank 2}

In the following Sections we consider the different Weyl equivalence
classes of braided vector spaces of unidentified type. We divide the
work depending on the dimension of such spaces. We consider in this
Section the 2-dimensional unidentified spaces, then some particular
cases in rank three and four, and finally the remaining cases, but
dividing the work in some families, according with the shape of the
associated generalized Dynkin diagram.

The unidentified braided vector spaces in \cite[Table 1]{H-classif
RS} of rank two are those in rows 7, 8, 9, 12, 13, 14, 15 and 16.

We will consider each possible row, describe the root system of each
braiding, and calculate the dimension of the corresponding Nichols
algebra.

\begin{remark}\label{remark:simple roots}
The hyperword associated to a simple root $\alpha_i$ is the one associated to the unique Lyndon word of degree $\alpha_i$:
$x_i$. Also, the hyperword associated to a root of the way $m\alpha_1+\alpha_2$ is $x_{m\alpha_1+\alpha_2}=(\ad_c x_1)^mx_2$. Moreover, for the braidings considered in this Subsection we have the following possible hyperwords:
\begin{align*}
& x_{\alpha_1+2\alpha_2}=\left[x_{\alpha_1+\alpha_2},x_2\right]_c,\quad & x_{3\alpha_1+2\alpha_2}=\left[x_{2\alpha_1+\alpha_2},x_{\alpha_1+\alpha_2}\right]_c,
\\ & x_{4\alpha_1+3\alpha_2}=\left[x_{3\alpha_1+2\alpha_2},x_{\alpha_1+\alpha_2}\right]_c, \quad & x_{5\alpha_1+2\alpha_2}=
\left[ x_{3\alpha_1+\alpha_2},x_{2\alpha_1+\alpha_2}\right]_c,
\\ & x_{5\alpha_1+3\alpha_2}=\left[x_{2\alpha_1+\alpha_2},x_{3\alpha_1+2\alpha_2}\right]_c,\quad & x_{5\alpha_1+4\alpha_2}
=\left[ x_{4\alpha_1+3\alpha_2},x_{\alpha_1+\alpha_2}\right]_c,
\\ & x_{7\alpha_1+2\alpha_2}=\left[x_{4\alpha_1+\alpha_2},x_{3\alpha_1+\alpha_2}\right]_c,\quad & x_{7\alpha_1+3\alpha_2}
=\left[ x_{5\alpha_1+2\alpha_2},x_{2\alpha_1+\alpha_2}\right]_c,
\\ & x_{7\alpha_1+4\alpha_2}=\left[x_{2\alpha_1+\alpha_2},x_{5\alpha_1+3\alpha_2}\right]_c,\quad & x_{8\alpha_1+3\alpha_2}
=\left[ x_{3\alpha_1+\alpha_2},x_{5\alpha_1+2\alpha_2}\right]_c,
\\ & x_{8\alpha_1+5\alpha_2}=\left[x_{5\alpha_1+3\alpha_2},x_{3\alpha_1+2\alpha_2}\right]_c.
\end{align*}
\end{remark}

\begin{remark}\label{remark:dimension}
If $V$, $W$ are two Weyl equivalent braided vector spaces of diagonal type, then $\dim \cB(V)=\dim \cB(W)$. It follows from the fact
$$ \Delta^{\schi}_+=s_p\left( \Delta^{\chi}\setminus\{\alpha_p\}\right) \cup \{\alpha_p\} $$
and that a hyperword of degree $\alpha$ has height $\ord \chi(\alpha,\alpha)$, so we calculate the dimension computing the number of terms of the PBW basis;
i.e. multiplying the orders of the associated scalars.
\end{remark}

\begin{example}{\textbf{Row 7.}}\label{example:rank 2,row 7}
These are the first unidentified braided vector spaces, for which
$\zeta \in \G'_{12}$ is primitive. The following diagram shows the
action of the Weyl groupoid, where we indicate the vertices $1,2$ by
$\circ$, $\bullet$, respectively. We omit those symmetries not
changing the Dynkin diagram.

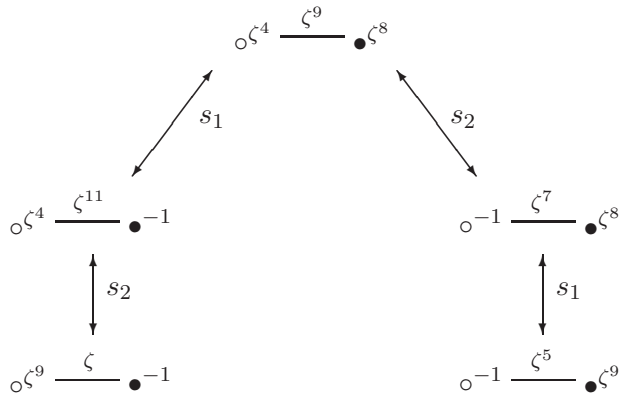
\begin{align*}
\begin{picture}(450,175)(10,10)
\put(180,150){\vector(-3,-4){30}} \put(150,110){\vector(3,4){30}}
\put(175,130){$\displaystyle s_1$} \put(280,110){\vector(-3,4){30}}
\put(250,150){\vector(3,-4){30}} \put(270,130){$\displaystyle s_2$}
\put(185,160){$\displaystyle\xymatrix{\circ^{\zeta^4}
\ar@{-}[r]^{\zeta^9} & \bullet^{\zeta^8} }$}
\put(100,90){$\displaystyle\xymatrix{\circ^{\zeta^4}
\ar@{-}[r]^{\zeta^{11}} & \bullet^{-1} }$}
\put(270,90){$\displaystyle\xymatrix{\circ^{-1} \ar@{-}[r]^{\zeta^7}
& \bullet^{\zeta^8} }$} \put(135,80){\vector(0,-1){30}}
\put(135,50){\vector(0,1){30}} \put(140,65){$\displaystyle s_2$}
\put(305,80){\vector(0,-1){30}} \put(305,50){\vector(0,1){30}}
\put(310,65){$\displaystyle s_1$}
\put(100,30){$\displaystyle\xymatrix{\circ^{\zeta^9}
\ar@{-}[r]^{\zeta} & \bullet^{-1} }$}
\put(270,30){$\displaystyle\xymatrix{\circ^{-1} \ar@{-}[r]^{\zeta^5}
& \bullet^{\zeta^9}}$}
\end{picture}
\end{align*}

\noindent In this case, $\cO(\chi)$ is empty. For each one of these
vector spaces $V$, $\dim\cB(V)=2^43^2=144.$

\noindent \vi $\xymatrix{\circ^{\zeta^4} \ar@{-}[r]^{\zeta^9} &
\bullet^{\zeta^8} }$. In this case,
$$ \Delta_+^{\chi}= \left\{ \alpha_1,2\alpha_1+\alpha_2,\alpha_1+\alpha_2,\alpha_1+2\alpha_2,\alpha_2 \right\}. $$
Considering the corresponding hyperwords from Remark \ref{remark:simple roots} and following Theorem \ref{thm:presentacion minima}, $\cB(V)$ has a
presentation by generators $x_1,x_2$, and relations
$$x_1^3=x_2^3= \left[x_1,x_{\alpha_1+2\alpha_2}\right]_c -\frac{\zeta^{10}(1-\zeta^7)q_{12}}{1-\zeta^9}x_{\alpha_1+\alpha_2}^2=0.$$

\noindent \vii $\xymatrix{\circ^{\zeta^4} \ar@{-}[r]^{\zeta^{11}} &
\bullet^{-1} }$. The set of positive roots is in this case
$$ \Delta_+^{\chi}= \left\{ \alpha_1,2\alpha_1+\alpha_2,3\alpha_1+2\alpha_2,\alpha_1+\alpha_2,\alpha_2 \right\}. $$
Therefore $\cB(V)$ has a presentation by generators $x_1,x_2$, and relations
$$ x_1^3=x_2^2=\left[x_{3\alpha_1+2\alpha_2},x_{12}\right]_c=0.$$

\noindent \viii $\xymatrix{\circ^{-1} \ar@{-}[r]^{\zeta^7} &
\bullet^{\zeta^8} }$. The positive roots are the same as in \vii
exchanging $1$ by $2$, and $\cB(V)$ admits an analogous
presentation.

\noindent \viv $\xymatrix{\circ^{\zeta^9} \ar@{-}[r]^{\zeta} &
\bullet^{-1} }$. The positive roots are
$$ \Delta_+^{\chi}= \left\{ \alpha_1,3\alpha_1+\alpha_2,2\alpha_1+\alpha_2,\alpha_1+\alpha_2,\alpha_2 \right\}. $$
Therefore $\cB(V)$ has a presentation by generators $x_1,x_2$, and relations
$$ x_1^4=x_2^2= \left[ x_{2\alpha_1+\alpha_2},x_{12}\right]_c=0.$$

\noindent \vv $\xymatrix{\circ^{-1} \ar@{-}[r]^{\zeta^5} &
\bullet^{\zeta^9} }$. The set of positive roots is the same as in
\viv exchanging $1$ by $2$, and $\cB(V)$ has an analogous
presentation.
\end{example}

\begin{example}{\textbf{Row 8.}}\label{example:rank 2,row 8}
Let $\zeta \in \G'_{12}$ be primitive. The roots in $\cO(\chi)$ are those such that $q_\alpha=\zeta^5$. For each one of these vector spaces $V$,
$\dim\cB(V)= 2^43^3=432.$

\noindent \vi $\xymatrix{\circ^{\zeta^8} \ar@{-}[r]^{\zeta} & \circ^{\zeta^8} }$. Their positive roots are
$$ \Delta_+^{\chi}= \left\{ \alpha_1,2\alpha_1+\alpha_2,\alpha_1+\alpha_2,\alpha_1+2\alpha_2,\alpha_2 \right\}. $$
According to the Theorem \ref{thm:presentacion minima}, $\cB(V)$ has a presentation by generators $x_1,x_2$, and relations
$$x_1^3=x_2^3=x_{\alpha_1+\alpha_2}^{12}=\left[x_1,x_{\alpha_1+2\alpha_2}\right]_c-\frac{\zeta^{10}(1-\zeta^7)q_{12}}{1-\zeta^9}x_{\alpha_1+\alpha_2}^2=0.$$

\noindent \vii $\xymatrix{\circ^{\zeta^8} \ar@{-}[r]^{\zeta^3} & \circ^{-1} }$. In this case,
$$ \Delta_+^{\chi}= \left\{ \alpha_1,2\alpha_1+\alpha_2,3\alpha_1+2\alpha_2,\alpha_1+\alpha_2,\alpha_2 \right\}. $$
With the corresponding hyperwords from Remark \ref{remark:simple roots}, $\cB(V)$ has a presentation by generators $x_1,x_2$, and relations
$$ x_1^3=x_2^2=x_{\alpha_1+\alpha_2}^{12}=\left[x_{3\alpha_1+2\alpha_2},x_{12}\right]_c=0.$$

\noindent \viii $\xymatrix{\circ^{\zeta^9} \ar@{-}[r]^{\zeta} & \circ^{-1} }$. The positive roots for this braiding are
$$ \Delta_+^{\chi}= \left\{ \alpha_1,3\alpha_1+\alpha_2,2\alpha_1+\alpha_2,\alpha_1+\alpha_2,\alpha_2 \right\}. $$
Therefore $\cB(V)$ has a presentation by generators $x_1,x_2$, and relations
$$ x_1^{12}=x_2^2=(\ad_cx_1)^4x_2= \left[ x_{2\alpha_1+\alpha_2},x_{12}\right]_c=0.$$
\end{example}

\begin{example}{\textbf{Row 9.}}\label{example:rank 2,row 9}
Let $\zeta \in \G'_9$ be primitive. The roots in $\cO(\chi)$ are those such that $N_\alpha=18$. For each one of these vector spaces $V$ we have $\dim\cB(V)=
2^43^6.$

\noindent \vi $\xymatrix{\circ^{-\zeta} \ar@{-}[r]^{\zeta^7} & \circ^{\zeta^3} }$. Their positive roots are
$$ \Delta_+^{\chi}= \left\{ \alpha_1,2\alpha_1+\alpha_2,3\alpha_1+2\alpha_2,\alpha_1+\alpha_2,\alpha_1+2\alpha_2,\alpha_2 \right\}. $$
Following Theorem \ref{thm:presentacion minima}, $\cB(V)$ has a presentation by generators $x_1,x_2$, and relations
$$x_1^{18}=x_2^3=x_{\alpha_1+\alpha_2}^{18}=\left[x_1,x_{\alpha_1+2\alpha_2}\right]_c+\frac{\zeta^5(1-\zeta)q_{12}}{1-\zeta^7}
x_{\alpha_1+\alpha_2}^2=0.$$

\noindent \vii $\xymatrix{\circ^{\zeta^3} \ar@{-}[r]^{\zeta^8} & \circ^{-1} }$. The set of positive roots is the following
$$ \Delta_+^{\chi}= \left\{ \alpha_1,2\alpha_1+\alpha_2,3\alpha_1+2\alpha_2,4\alpha_1+3\alpha_2,\alpha_1+\alpha_2,\alpha_2 \right\}. $$
Therefore $\cB(V)$ has a presentation by generators $x_1,x_2$, and relations
$$ x_1^3=x_2^2=x_{\alpha_1+\alpha_2}^{18}=x_{2\alpha_1+\alpha_2}^{18}=\left[x_{2\alpha_1+\alpha_2},x_{3\alpha_1+2\alpha_2}\right]_c=0.$$

\noindent \viii $\xymatrix{\circ^{-\zeta^2} \ar@{-}[r]^{\zeta} & \circ^{-1} }$. In this case,
$$ \Delta_+^{\chi}= \left\{\alpha_1,4\alpha_1+\alpha_2,3\alpha_1+\alpha_2,2\alpha_1+\alpha_2,\alpha_1+\alpha_2,\alpha_2 \right\}. $$
Therefore $\cB(V)$ has a presentation by generators $x_1,x_2$, and relations
$$ x_1^{18}=x_2^2=x_{2\alpha_1+\alpha_2}^{18}=(\ad_cx_1)^5x_2= \left[x_{2\alpha_1+\alpha_2},x_{12}\right]_c=0.$$
\end{example}

\begin{example}{\textbf{Row 12.}}\label{example:rank 2,row 12}
Let $\zeta \in \G'_{24}$ be primitive. Notice that the roots in $\cO(\chi)$ are those such that $N_\alpha=24$. For each one of these vector spaces $V$, we
have $ \dim \cB(V)= 2^{10}3^4$.

\noindent \vi $\xymatrix{\circ^{\zeta^6}\ar@{-}[r]^{\zeta^{11}}&\circ^{\zeta^8} }$. The set of positive roots is
$$ \left\{\alpha_1,3\alpha_1+\alpha_2,2\alpha_1+\alpha_2,3\alpha_1+2\alpha_2,4\alpha_1+3\alpha_2,\alpha_1+\alpha_2,\alpha_1+2\alpha_2,
\alpha_2 \right\}. $$
Considering the associated hyperwords, Theorem \ref{thm:presentacion minima} establishes that $\cB(V)$ has a presentation by generators $x_1,x_2$, and
relations
\begin{align*}
& x_1^4=x_2^3= x_{3\alpha_1+\alpha_2}^{24}=x_{\alpha_1+\alpha_2}^{24}=0,
\\ & (1-\zeta^{11})\zeta^4q_{21} \left[x_1,x_{\alpha_1+2\alpha_2}\right]_c = (1-\zeta^{19}) x_{\alpha_1+\alpha_2}^2.
\end{align*}

\noindent \vii $\xymatrix{\circ^{\zeta^6} \ar@{-}[r]^{\zeta} & \circ^{\zeta^{-1}} }$. The set of positive roots is in this case:
$$ \left\{\alpha_1,3\alpha_1+\alpha_2,5\alpha_1+2\alpha_2,2\alpha_1+\alpha_2,5\alpha_1+3\alpha_2,3\alpha_1+2\alpha_2,\alpha_1+\alpha_2,
\alpha_2 \right\}. $$
Following Theorem \ref{thm:presentacion minima}, $\cB(V)$ has a presentation by generators $x_1,x_2$, and relations
$$ x_1^4=x_2^{24}= x_{2\alpha_1+\alpha_2}^{24}=(\ad_c x_2)^2x_1=\left[x_{3\alpha_1+2\alpha_2},x_{12}\right]_c =0.$$

\noindent \viii $\xymatrix{\circ^{\zeta^8} \ar@{-}[r]^{\zeta^5} & \circ^{-1} }$. The positive roots are in this case
$$ \left\{\alpha_1,2\alpha_1+\alpha_2,5\alpha_1+3\alpha_2,3\alpha_1+2\alpha_2,4\alpha_1+3\alpha_2,5\alpha_1+4\alpha_2,\alpha_1+\alpha_2,
\alpha_2 \right\}. $$
Therefore $\cB(V)$ has a presentation by generators $x_1,x_2$, and relations
\begin{align*}
& x_1^3=x_2^2=x_{5\alpha_1+3\alpha_2}^{24}=x_{\alpha_1+\alpha_2}^{24}=0,
\\ & \left[x_{2\alpha_1+\alpha_2},x_{4\alpha_1+3\alpha_2}\right]_c = \frac{1+\zeta+\zeta^6+2\zeta^7+\zeta^{17}}{(1+\zeta^4+\zeta^6+\zeta^{11})\zeta^{10}q_{21}}
\ x_{3\alpha_1+2\alpha_2}^2.
\end{align*}

\noindent \viv $\xymatrix{\circ^{\zeta} \ar@{-}[r]^{\zeta^{19}} & \circ^{-1} }$. The set of positive roots is in this case:
$$ \left\{\alpha_1,5\alpha_1+\alpha_2,4\alpha_1+\alpha_2,3\alpha_1+\alpha_2,5\alpha_1+2\alpha_2,2\alpha_1+\alpha_2,\alpha_1+\alpha_2,
\alpha_2 \right\}. $$
Therefore $\cB(V)$ has a presentation by generators $x_1,x_2$, and relations
$$x_1^{24}=x_2^2= x_{5\alpha_1+2\alpha_2}^{24}=(\ad_c x_1)^6x_2=\left[ x_{2\alpha_1+\alpha_2} , x_{12}\right]_c=0.$$
\end{example}

\begin{example}{\textbf{Row 13.}}\label{example:rank 2,row 13}
Let $\zeta \in \G'_5$ be primitive. The roots in $\cO(\chi)$ are those such that $q_\alpha\in\{\zeta,-\zeta^3\}$. For each one of these vector spaces $V$,
we have $ \dim \cB(V)= 2^65^4$.

\noindent \vi $\xymatrix{\circ^{\zeta}\ar@{-}[r]^{\zeta^2}&\circ^{-1} }$.  The set of positive roots is in this case:
$$ \left\{\alpha_1,3\alpha_1+\alpha_2,2\alpha_1+\alpha_2,5\alpha_1+3\alpha_2,3\alpha_1+2\alpha_2,4\alpha_1+3\alpha_2,\alpha_1+\alpha_2,
\alpha_2 \right\}. $$
Following Theorem \ref{thm:presentacion minima}, $\cB(V)$ has a presentation by generators $x_1,x_2$, and relations
$$ x_1^5=x_2^2= x_{2\alpha_1+\alpha_2}^{10}=x_{3\alpha_1+2\alpha_2}^5=x_{\alpha_1+\alpha_2}^{10}=(\ad_cx_1)^4x_2=\left[x_{4\alpha_1+3\alpha_2},x_{12}\right]_c=0.$$

\noindent \vii $\xymatrix{\circ^{-\zeta^3} \ar@{-}[r]^{\zeta^3} & \circ^{-1} }$.  For this braiding, the positive roots are
$$ \left\{\alpha_1,4\alpha_1+\alpha_2,3\alpha_1+\alpha_2,5\alpha_1+2\alpha_2,2\alpha_1+\alpha_2,3\alpha_1+2\alpha_2,\alpha_1+\alpha_2,
\alpha_2 \right\}. $$
Therefore $\cB(V)$ has a presentation by generators $x_1,x_2$, and relations
\begin{align*}
& x_1^{10}=x_2^2= x_{3\alpha_1+\alpha_2}^5=x_{\alpha_1+\alpha_2}^5=x_{2\alpha_1+\alpha_2}^{10}=0,
\\ & (\ad_cx_1)^5x_2=\left[x_1,x_{3\alpha_1+2\alpha_2}\right]_c+q_{12}x_{2\alpha_1+\alpha_2}^2=0.
\end{align*}
\end{example}

\begin{example}{\textbf{Row 14.}}\label{example:rank 2,row 14}
Let $\zeta \in \G'_{20}$ be primitive. The roots in $\cO(\chi)$ are those such that $N_\alpha=20$. For each one of these vector spaces $V$,
we have $ \dim \cB(V)= 2^85^4$.

\noindent \vi $\xymatrix{\circ^{\zeta}\ar@{-}[r]^{\zeta^{17}}&\circ^{-1} }$.  The set of positive roots is the same as in Example
\ref{example:rank 2,row 13},\vi.
Following Theorem \ref{thm:presentacion minima}, $\cB(V)$ has a presentation by generators $x_1,x_2$, and relations
$$ x_1^{20}=x_2^2= x_{3\alpha_1+2\alpha_2}^{20}=(\ad_cx_1)^4x_2=\left[x_{4\alpha_1+3\alpha_2},x_{12}\right]_c=0.$$

\noindent \vii $\xymatrix{\circ^{\zeta^{11}} \ar@{-}[r]^{\zeta^7} & \circ^{-1} }$. The set of positive roots is again the same as in Example
\ref{example:rank 2,row 13},\vi.
According to Theorem \ref{thm:presentacion minima}, $\cB(V)$ has a presentation by generators $x_1,x_2$, and relations
$$ x_1^{20}=x_2^2= x_{3\alpha_1+2\alpha_2}^{20}=(\ad_cx_1)^4x_2=\left[x_{4\alpha_1+3\alpha_2},x_{12}\right]_c=0.$$

\noindent \viii $\xymatrix{\circ^{\zeta^8}\ar@{-}[r]^{\zeta^3}&\circ^{-1} }$. The set of positive roots is the same as in Example
\ref{example:rank 2,row 13},\vii.
Following Theorem \ref{thm:presentacion minima}, $\cB(V)$ has a presentation by generators $x_1,x_2$, and relations
$$  x_1^5=x_2^2=x_{3\alpha_1+\alpha_2}^{20}=x_{\alpha_1+\alpha_2}^{20}=
\left[x_1,x_{3\alpha_1+2\alpha_2}\right]_c+\frac{(1-\zeta^{17})q_{12}}{1-\zeta^2}x_{2\alpha_1+\alpha_2}^2=0.$$

\noindent \viv $\xymatrix{\circ^{\zeta^8} \ar@{-}[r]^{\zeta^{13}} & \circ^{-1} }$. The set of positive roots is again the same as in Example
\ref{example:rank 2,row 13},\vii.
According to Theorem \ref{thm:presentacion minima}, $\cB(V)$ has a presentation by generators $x_1,x_2$, and relations
$$ x_1^5=x_2^2=x_{3\alpha_1+\alpha_2}^{20}=x_{\alpha_1+\alpha_2}^{20}=
\left[x_1,x_{3\alpha_1+2\alpha_2}\right]_c+\frac{(1-\zeta^7)q_{12}}{1-\zeta^2}x_{2\alpha_1+\alpha_2}^2=0.$$
\end{example}

\begin{example}{\textbf{Row 15.}}\label{example:rank 2,row 15}
Let $\zeta \in \G'_{15}$ be primitive. Notice that the roots in $\cO(\chi)$ are those such that $N_\alpha=30$. For each one of these vector spaces $V$, we
have $\dim \cB(V)= 2^43^45^4=30^4$.

\noindent \vi $\xymatrix{\circ^{-\zeta}\ar@{-}[r]^{-\zeta^{12}}&\circ^{\zeta^5} }$. The set of positive roots is
$$ \left\{\alpha_1,3\alpha_1+\alpha_2,5\alpha_1+2\alpha_2,2\alpha_1+\alpha_2,3\alpha_1+2\alpha_2,\alpha_1+\alpha_2,\alpha_1+2\alpha_2,
\alpha_2 \right\}. $$
Considering the associated hyperwords, Theorem \ref{thm:presentacion minima} says that $\cB(V)$ has a presentation by generators $x_1,x_2$, and
relations
\begin{align*}
& x_1^{30}=x_2^3=x_{3\alpha_1+\alpha_2}^{30}=(\ad_cx_1)^4x_2=0,
\\ & \left[x_1,x_{\alpha_1+2\alpha_2}\right]_c +\frac{\zeta^{10}(1+\zeta^{13})q_{12}}{1+\zeta^{12}} x_{\alpha_1+\alpha_2}^2=
\left[x_{3\alpha_1+2\alpha_2},x_{12}\right]_c =0.
\end{align*}

\noindent \vii $\xymatrix{\circ^{\zeta^3} \ar@{-}[r]^{-\zeta^4} & \circ^{-\zeta^{11}} }$. The set of positive roots is in this case:
$$ \left\{\alpha_1,4\alpha_1+\alpha_2,3\alpha_1+\alpha_2,2\alpha_1+\alpha_2,3\alpha_1+2\alpha_2,4\alpha_1+3\alpha_2,\alpha_1+\alpha_2,
\alpha_2 \right\}. $$
Following Theorem \ref{thm:presentacion minima}, $\cB(V)$ has a presentation by generators $x_1,x_2$, and relations
\begin{align*}
& x_1^5=x_2^{30}=x_{2\alpha_1+\alpha_2}^{30}=(\ad_cx_2)^2x_1=0,
\\ & \left[x_1,x_{3\alpha_1+2\alpha_2}\right]_c -\frac{(1-\zeta^2)\zeta^9q_{12}}{1+\zeta^7} x_{2\alpha_1+\alpha_2}^2=
\left[x_{4\alpha_1+3\alpha_2},x_{12}\right]_c =0.
\end{align*}

\noindent \viii $\xymatrix{\circ^{\zeta^5} \ar@{-}[r]^{-\zeta^{13}} & \circ^{-1} }$. The positive roots are in this case
$$ \left\{\alpha_1,2\alpha_1+\alpha_2,5\alpha_1+3\alpha_2,8\alpha_1+5\alpha_2,3\alpha_1+2\alpha_2,4\alpha_1+3\alpha_2,\alpha_1+\alpha_2,
\alpha_2 \right\}. $$
Therefore $\cB(V)$ has a presentation by generators $x_1,x_2$, and relations
$$ x_1^3=x_2^2=x_{2\alpha_1+\alpha_2}^{30}=x_{4\alpha_1+3\alpha_2}^{30}=\left[x_{4\alpha_1+3\alpha_2},x_{12}\right]_c =0.$$

\noindent \viv $\xymatrix{\circ^{\zeta^3} \ar@{-}[r]^{-\zeta^2} & \circ^{-1} }$. The set of positive roots is in this case:
$$ \left\{\alpha_1,4\alpha_1+\alpha_2,3\alpha_1+\alpha_2,8\alpha_1+3\alpha_2,5\alpha_1+2\alpha_2,2\alpha_1+\alpha_2,\alpha_1+\alpha_2,
\alpha_2 \right\}. $$
Therefore $\cB(V)$ has a presentation by generators $x_1,x_2$, and relations
$$ x_1^5=x_2^2=x_{4\alpha_1+\alpha_2}^{30}=x_{2\alpha_1+\alpha_2}^{30}=\left[x_{2\alpha_1+\alpha_2},x_{12}\right]_c
=\left[x_{5\alpha_1+2\alpha_2},x_{112}\right]_c =0.$$
\end{example}

\begin{example}{\textbf{Row 16.}}\label{example:rank 2,row 16}
Let $\zeta \in \G'_7$ be primitive. The roots in $\cO(\chi)$ are those such that $q_\alpha\in\{-\zeta,-\zeta^5\}$. For each one of these vector spaces $V$,
we have $ \dim \cB(V)= 2^{12}7^6$.

\noindent \vi $\xymatrix{\circ^{\zeta}\ar@{-}[r]^{\zeta^2}&\circ^{-1} }$.  The set of positive roots is in this case:
\begin{align*}
\Delta_+^{\chi}=& \left\{\alpha_1,3\alpha_1+\alpha_2,2\alpha_1+\alpha_2,7\alpha_1+4\alpha_2,5\alpha_1+3\alpha_2,8\alpha_1+5\alpha_2, \right.
\\ & \left. 3\alpha_1+2\alpha_2,7\alpha_1+5\alpha_2,4\alpha_1+3\alpha_2,5\alpha_1+4\alpha_2,\alpha_1+\alpha_2,\alpha_2 \right\}.
\end{align*}
Following Theorem \ref{thm:presentacion minima}, $\cB(V)$ has a presentation by generators $x_1,x_2$, and relations
\begin{align*}
& x_2^2= x_{\alpha}^{14}=0, \qquad \alpha= \alpha_1,2\alpha_1+\alpha_2,5\alpha_1+3\alpha_2,3\alpha_1+2\alpha_2,4\alpha_1+3\alpha_2,
\alpha_1+\alpha_2,
\\ & (\ad_cx_1)^4x_2=[x_{2\alpha_i+\alpha_j},x_{4\alpha_i+3\alpha_j}]_c-\frac{(1-\zeta+3\zeta^4-\zeta^6)q_{ij}}
{2\zeta+\zeta^2-\zeta^3-\zeta^4+\zeta^5-2} x_{3\alpha_i+2\alpha_j}^2 =0,
\end{align*}

\noindent \vii $\xymatrix{\circ^{-\zeta^3} \ar@{-}[r]^{\zeta^3} & \circ^{-1} }$.  For this braiding, the positive roots are
\begin{align*}
\Delta_+^{\chi}=& \left\{\alpha_1,5\alpha_1+\alpha_2,4\alpha_1+\alpha_2,7\alpha_1+2\alpha_2,3\alpha_1+\alpha_2,8\alpha_1+3\alpha_2, \right.
\\ & \left. 5\alpha_1+2\alpha_2,7\alpha_1+3\alpha_2,2\alpha_1+\alpha_2,3\alpha_1+2\alpha_2,\alpha_1+\alpha_2,\alpha_2 \right\}.
\end{align*}
Therefore $\cB(V)$ has a presentation by generators $x_1,x_2$, and relations
\begin{align*}
& x_2^2= x_{\alpha}^{14}=0, \qquad \alpha= \alpha_1,4\alpha_1+\alpha_2,3\alpha_1+\alpha_2,5\alpha_1+2\alpha_2,2\alpha_1+\alpha_2,\alpha_1+\alpha_2,
\\ & (\ad_cx_1)^6x_2= \left[x_1,x_{3\alpha_1+2\alpha_2}\right]_c +\frac{(1+\zeta^4)q_{12}}{1-\zeta^2}x_{2\alpha_1+\alpha_2}^2=0,
\end{align*}
\end{example}

\medskip

\section{Examples in rank 3 and 4}

Now we consider three Weyl equivalences classes in rank three and
two in rank four, and make the same work as in the previous Section.

\begin{example}{\textbf{Rank 3, row 13.}}\label{example:rank 3,row 13}
In this case we have two different Weyl groupoids, depending on the
order of $\zeta$, with the same root systems (the associated
braidings are different). We analize each case.
\smallskip

\noindent $\mathbf{1.}$ Let $\zeta \in \G'_3$ be primitive. The roots
in $\cO(\chi)$ are those such that $\ord q_\alpha\in\{3,6\}$. For
each one of these vector spaces $V$, $\dim\cB(V)= 2^4 3^6 6^3=2^7
3^9$.

\noindent \vi
$\xymatrix{\circ^{\zeta}\ar@{-}[r]^{\zeta^2}&\circ^{\zeta}\ar@{-}[r]^{\zeta}&\circ^{-1}}$.
The set of positive roots is the following:
\begin{align*}
\Delta_+^{\chi}=&
\left\{\alpha_1,\alpha_2,\alpha_3,\alpha_1+\alpha_2,\alpha_2+\alpha_3,
2\alpha_1+2\alpha_2+\alpha_3,\alpha_1+\alpha_2+\alpha_3,\alpha_1+2\alpha_2+\alpha_3,\right.
\\ & \left. \alpha_1+2\alpha_2+2\alpha_3,2\alpha_2+\alpha_3,\alpha_1+3\alpha_2+2\alpha_3,
2\alpha_1+3\alpha_2+2\alpha_3,2\alpha_1+4\alpha_2+3\alpha_3\right\}.
\end{align*}
According to Theorem \ref{thm:presentacion minima}, $\cB(V)$ has a
presentation by generators $x_1,x_2,x_3$, and relations
\begin{align*}
& x_3^2= (\ad_cx_1)^2x_2= (\ad_cx_1)x_3= (\ad_cx_2)^2x_1=
[x_{223},x_{23}]_c = \left[[x_{123},x_2]_c,x_2\right]_c =0,
\\ &x_{\alpha}^{3}=0, \qquad \alpha= \alpha_1,\alpha_1+\alpha_2,\alpha_2,\alpha_1+2\alpha_2+2\alpha_3,\alpha_1+3\alpha_2+2\alpha_3,
2\alpha_1+3\alpha_2+2\alpha_3,
\\ & x_{\alpha}^{6}=0, \qquad \alpha=
\alpha_2+\alpha_3,\alpha_1+\alpha_2+\alpha_3,\alpha_1+2\alpha_2+\alpha_3.
\end{align*}

\noindent \vii
$\xymatrix{\circ^{\zeta}\ar@{-}[r]^{\zeta^2}&\circ^{-\zeta^2}\ar@{-}[r]^{\zeta^2}&\circ^{-1}}$.
In this case, we have
\begin{align*}
\Delta_+^{\chi}=&
\left\{\alpha_1,\alpha_2,\alpha_3,\alpha_1+\alpha_2,\alpha_1+2\alpha_2,
2\alpha_1+2\alpha_2+\alpha_3,\alpha_1+\alpha_2+\alpha_3,2\alpha_1+4\alpha_2+\alpha_3,\right.
\\ & \left. 2\alpha_2+\alpha_3,\alpha_2+\alpha_3,\alpha_1+3\alpha_2+\alpha_3,
2\alpha_1+3\alpha_2+\alpha_3,\alpha_1+2\alpha_2+\alpha_3\right\}.
\end{align*}
Theorem \ref{thm:presentacion minima} says that $\cB(V)$ is
presented by generators $x_1,x_2,x_3$, and relations
\begin{align*}
& x_3^2= (\ad_cx_1)^2x_2= (\ad_cx_1)x_3= (\ad_cx_2)^3x_1=
(\ad_cx_2)^3x_3=0,
\\ & [x_3,x_{221}]_c
+q_{21}[x_{321},x_2]_c+q_{32}\zeta^2(1-\zeta^2)x_2x_{321} =0,
\\ &x_{\alpha}^{3}=0, \qquad \alpha= \alpha_1,\alpha_1+\alpha_2+\alpha_3,\alpha_2+\alpha_3,\alpha_1+2\alpha_2,
\alpha_1+3\alpha_2+\alpha_3,2\alpha_1+3\alpha_2+\alpha_3,
\\ & x_{\alpha}^{6}=0, \qquad \alpha=
\alpha_2,\alpha_1+\alpha_2,\alpha_1+2\alpha_2+\alpha_3.
\end{align*}
\smallskip

\noindent $\mathbf{2.}$ Now, let $\zeta \in \G'_6$ be primitive. The
roots in $\cO(\chi)$ are those such that $\ord q_\alpha=6$. In this
case, $ \dim \cB(V)= 2^4 3^3 6^6=2^{10} 3^9$.

\noindent \vi
$\xymatrix{\circ^{\zeta}\ar@{-}[r]^{\zeta^5}&\circ^{\zeta}\ar@{-}[r]^{\zeta^4}&\circ^{-1}}$.
The set of positive roots is the same as in $\mathbf{1.}$\vi. Using
Theorem \ref{thm:presentacion minima}, we deduce that $\cB(V)$ has a
presentation by generators $x_1,x_2,x_3$, and relations
\begin{align*}
& x_3^2= (\ad_cx_1)^2x_2= (\ad_cx_1)x_3= (\ad_cx_2)^2x_1=
(\ad_cx_2)^3x_2=0,
\\ &x_{\alpha}^{6}=0, \qquad \alpha= \alpha_1,\alpha_1+\alpha_2,\alpha_2,\alpha_1+2\alpha_2+2\alpha_3,\alpha_1+3\alpha_2+2\alpha_3,
2\alpha_1+3\alpha_2+2\alpha_3.
\end{align*}

\noindent \vii
$\xymatrix{\circ^{\zeta}\ar@{-}[r]^{\zeta^2}&\circ^{-\zeta^2}\ar@{-}[r]^{\zeta^2}&\circ^{-1}}$.
The positive roots are the same as in $\mathbf{1.}$\vii. By Theorem
\ref{thm:presentacion minima}, $\cB(V)$ is presented by generators
$x_1,x_2,x_3$, and relations
\begin{align*}
& x_3^2=x_2^3=(\ad_cx_1)^2x_2 = (\ad_cx_1)x_3= [x_{223},x_{23}]_c=0,
\\ & [x_1,x_{223}]_c
+q_{23}[x_{123},x_2]_c-q_{12}x_2x_{123} =0,
\\ &x_{\alpha}^{6}=0, \qquad \alpha= \alpha_1,\alpha_1+\alpha_2+\alpha_3,\alpha_2+\alpha_3,\alpha_1+2\alpha_2,
\alpha_1+3\alpha_2+\alpha_3,2\alpha_1+3\alpha_2+\alpha_3.
\end{align*}
\end{example}

\begin{example}{\textbf{Rank 3, row 17.}}\label{example:rank 3,row 17}
Let $\zeta\in\G'_3$ be primitive. The following diagram shows the
action of the Weyl groupoid, where we indicate the vertices $1,2,3$
by $\circ$, $\ast$, $\bullet$, respectively, and we omit those
symmetries which do not change the Dynkin diagram.

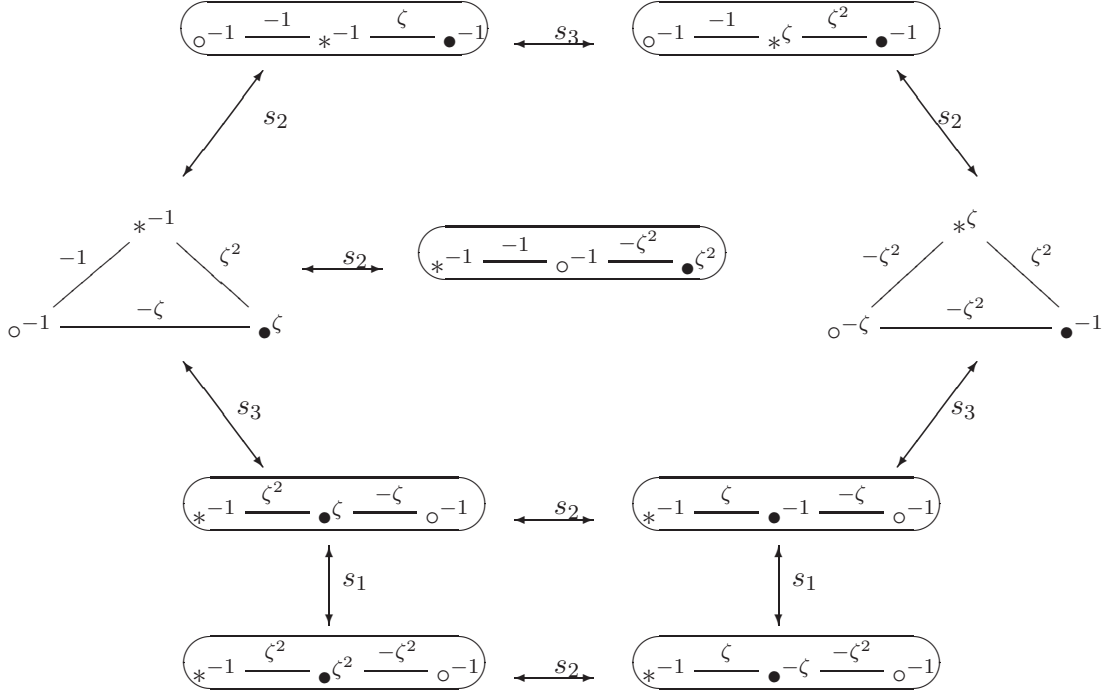
\begin{align*}
\begin{picture}(450,300)(10,10)
\put(80,270){$\displaystyle\xymatrix{\circ^{-1}\ar@{-}[r]^{-1}&\ast^{-1}
\ar@{-}[r]^{\zeta}&\bullet^{-1}}$}
\put(137,276){\oval(116,20)}
\put(205,270){\vector(1,0){30}}
\put(235,270){\vector(-1,0){30}}
\put(220,272){$\displaystyle s_3$}
\put(250,270){$\displaystyle\xymatrix{\circ^{-1}\ar@{-}[r]^{-1}&\ast^{\zeta}
\ar@{-}[r]^{\zeta^2}&\bullet^{-1}}$}
\put(308,276){\oval(116,20)}
\put(110,260){\vector(-3,-4){30}}
\put(80,220){\vector(3,4){30}}
\put(110,240){$\displaystyle s_2$}
\put(380,220){\vector(-3,4){30}}
\put(350,260){\vector(3,-4){30}}
\put(365,240){$\displaystyle s_2$}
\put(170,185){$\displaystyle\xymatrix{\ast^{-1}\ar@{-}[r]^{-1}&\circ^{-1}
\ar@{-}[r]^{-\zeta^2}&\bullet^{\zeta^2}}$}
\put(227,191){\oval(116,20)}
\put(125,185){\vector(1,0){30}}
\put(155,185){\vector(-1,0){30}}
\put(140,187){$\displaystyle s_2$}
\put(10,200){$\displaystyle\xymatrix{ & \ast^{-1}
\ar@{-}[rd]^{\zeta^2}&\\\circ^{-1}
\ar@{-}[ru]^{-1}\ar@{-}[rr]^{-\zeta} & & \bullet^{\zeta} }$}
\put(80,150){\vector(3,-4){30}}
\put(110,110){\vector(-3,4){30}}
\put(100,130){$\displaystyle s_3$}
\put(320,200){$\displaystyle\xymatrix{ & \ast^{\zeta}
\ar@{-}[rd]^{\zeta^2}&\\\circ^{-\zeta}
\ar@{-}[ru]^{-\zeta^2}\ar@{-}[rr]^{-\zeta^2} &  & \bullet^{-1} }$}
\put(380,150){\vector(-3,-4){30}}
\put(350,110){\vector(3,4){30}}
\put(370,130){$\displaystyle s_3$}
\put(80,90){$\displaystyle\xymatrix{\ast^{-1}\ar@{-}[r]^{\zeta^2}&\bullet^{\zeta}
\ar@{-}[r]^{-\zeta}&\circ^{-1}}$}
\put(137,96){\oval(114,20)}
\put(205,90){\vector(1,0){30}}
\put(235,90){\vector(-1,0){30}}
\put(220,92){$\displaystyle s_2$}
\put(250,90){$\displaystyle\xymatrix{\ast^{-1}\ar@{-}[r]^{\zeta}&\bullet^{-1}
\ar@{-}[r]^{-\zeta}&\circ^{-1}}$}
\put(308,96){\oval(116,20)}
\put(135,80){\vector(0,-1){30}}
\put(135,50){\vector(0,1){30}}
\put(140,65){$\displaystyle s_1$}
\put(305,80){\vector(0,-1){30}}
\put(305,50){\vector(0,1){30}}
\put(310,65){$\displaystyle s_1$}
\put(80,30){$\displaystyle\xymatrix{\ast^{-1}\ar@{-}[r]^{\zeta^2}&\bullet^{\zeta^2}
\ar@{-}[r]^{-\zeta^2}&\circ^{-1}}$}
\put(137,36){\oval(114,20)}
\put(205,30){\vector(1,0){30}}
\put(235,30){\vector(-1,0){30}}
\put(220,32){$\displaystyle s_2$}
\put(250,30){$\displaystyle\xymatrix{\ast^{-1}\ar@{-}[r]^{\zeta}&\bullet^{-\zeta}
\ar@{-}[r]^{-\zeta^2}&\circ^{-1}}$}
\put(308,36){\oval(116,20)}
\end{picture}
\end{align*}
In this case, the roots in $\cO(\chi)$ are those such that $N_\alpha=6$, and $ \dim \cB(V)=2^73^36=2^83^4$.

\noindent \vi $\xymatrix{\circ^{-1}\ar@{-}[r]^{-1}&\ast^{-1}\ar@{-}[r]^{\zeta}&\bullet^{-1}}$. We have that
\begin{align*}
\Delta_+^{\chi}=& \left\{\alpha_1,\alpha_2,\alpha_3,\alpha_1+\alpha_2,\alpha_2+\alpha_3,\alpha_1+\alpha_2+\alpha_3,
\alpha_1+2\alpha_2+\alpha_3,\alpha_1+2\alpha_2+2\alpha_3, \right.
\\& \left. \alpha_1+3\alpha_2+2\alpha_3,2\alpha_1+3\alpha_2+2\alpha_3,2\alpha_1+4\alpha_2+3\alpha_3\right\}.
\end{align*}
Following Theorem \ref{thm:presentacion minima}, $\cB(V)$ has a presentation by generators $x_1$, $x_2$, $x_3$, and relations
$$ x_1^2=x_2^2=x_3^2= x_{12}^2=(\ad_cx_1)x_3=\left[x_{32}, \left[x_{32},x_{321} \right]_c \right]_c = x_{\alpha_1+2\alpha_2+2\alpha_3}^6=0.$$

\noindent \vii $\xymatrix{\circ^{-1}\ar@{-}[r]^{-1}&\ast^{\zeta}\ar@{-}[r]^{\zeta^2}&\bullet^{-1}}$. For this braiding,
\begin{align*}
\Delta_+^{\chi}=& \left\{\alpha_1,\alpha_2,\alpha_3,\alpha_1+\alpha_2,\alpha_1+2\alpha_2,\alpha_2+\alpha_3,
\alpha_1+\alpha_2+\alpha_3,\alpha_1+3\alpha_2+\alpha_3, \right.
\\& \left. 2\alpha_1+3\alpha_2+\alpha_3,2\alpha_1+4\alpha_2+\alpha_3,\alpha_1+2\alpha_2+\alpha_3\right\}.
\end{align*}
By Theorem \ref{thm:presentacion minima}, $\cB(V)$ has a presentation by generators $x_1$, $x_2$, $x_3$, and relations
$$ x_1^2=x_2^3=x_3^2= (\ad_cx_2)^2x_3 =(\ad_cx_1)x_3=\left[x_{221}, x_{21}\right]_c =\left[x_{12},x_{123}\right]_c=x_{\alpha_1+2\alpha_2}^6=0.$$

\noindent \viii  $\xymatrix{ & \ast^{\zeta} \ar@{-}[rd]^{\zeta^2}&\\\circ^{-\zeta} \ar@{-}[ru]^{-\zeta^2}\ar@{-}[rr]^{-\zeta^2} &  & \bullet^{-1} }$. The root system is
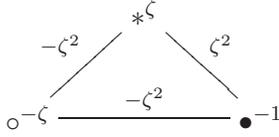
\begin{align*}
\Delta_+^{\chi}=& \left\{\alpha_1,\alpha_2,\alpha_3,\alpha_1+2\alpha_2,\alpha_1+\alpha_3,\alpha_1+\alpha_2,
\alpha_2+\alpha_3,2\alpha_1+\alpha_2+\alpha_3, \right.
\\& \left. \alpha_1+2\alpha_2+\alpha_3,\alpha_1+\alpha_2+\alpha_3,2\alpha_1+2\alpha_2+\alpha_3\right\}.
\end{align*}
By Theorem \ref{thm:presentacion minima}, $\cB(V)$ admits a presentation by generators $x_1$, $x_2$, $x_3$, and relations
\begin{align*}
x_1^6&=x_2^3=x_3^2= (\ad_cx_2)^2x_3 =(\ad_cx_1)^2x_2=(\ad_cx_1)^2x_3= 0
\\ x_{123}&=-q_{23}(1-\zeta^2)[x_{13},x_2]_c+q_{12}(1-\zeta^2)x_2x_{13}.
\end{align*}

\noindent \viv $\xymatrix{\ast^{-1}\ar@{-}[r]^{\zeta}&\bullet^{-1}\ar@{-}[r]^{-\zeta}&\circ^{-1}}$. In this case,
\begin{align*}
\Delta_+^{\chi}=& \left\{\alpha_1,\alpha_2,\alpha_3,\alpha_2+\alpha_3,\alpha_3+\alpha_1,\alpha_1+2\alpha_2+2\alpha_3,
\alpha_1+\alpha_2+\alpha_3,2\alpha_2+3\alpha_3+\alpha_1, \right.
\\& \left. \alpha_2+2\alpha_3+\alpha_1,\alpha_2+2\alpha_3+2\alpha_1,2\alpha_2+3\alpha_3+2\alpha_1\right\}.
\end{align*}
Following Theorem \ref{thm:presentacion minima}, $\cB(V)$ has a presentation by generators $x_1$, $x_2$, $x_3$, and relations
$$ x_1^2=x_2^2=x_3^2= (\ad_cx_1)x_2=\left[ \left[x_{13},x_{132} \right]_c,x_3\right]_c= x_{\alpha_1+\alpha_3}^6=0. $$

\noindent \vv $\xymatrix{\ast^{-1}\ar@{-}[r]^{\zeta}&\bullet^{-\zeta}\ar@{-}[r]^{-\zeta^2}&\circ^{-1}}$. Its root system is
\begin{align*}
\Delta_+^{\chi}=& \left\{\alpha_1,\alpha_2,\alpha_3,\alpha_2+\alpha_3,\alpha_2+2\alpha_3,\alpha_3+\alpha_1,
2\alpha_2+2\alpha_3+\alpha_1,\alpha_2+\alpha_3+\alpha_1, \right.
\\& \left. 2\alpha_2+3\alpha_3+\alpha_1,\alpha_2+2\alpha_3+\alpha_1,2\alpha_2+3\alpha_3+2\alpha_1\right\}.
\end{align*}
Following Theorem \ref{thm:presentacion minima}, $\cB(V)$ has a presentation by generators $x_1$, $x_2$, $x_3$, and relations
$$ x_1^2=x_2^2=x_3^6= (\ad_cx_3)^3x_2 =(\ad_cx_3)^2x_1=(\ad_cx_1)x_2= 0.$$

\noindent \vvi $\xymatrix{\ast^{-1}\ar@{-}[r]^{\zeta^2}&\bullet^{\zeta^2}\ar@{-}[r]^{-\zeta^2}&\circ^{-1}}$. The corresponding root system is
\begin{align*}
\Delta_+^{\chi}=& \left\{\alpha_1,\alpha_2,\alpha_3,\alpha_2+\alpha_3,\alpha_2+2\alpha_3,2\alpha_3+\alpha_1,
\alpha_3+\alpha_1,\alpha_2+\alpha_3+\alpha_1, \right.
\\& \left. \alpha_2+3\alpha_3+\alpha_1,\alpha_2+3\alpha_3+2\alpha_1,\alpha_2+2\alpha_3+\alpha_1\right\}.
\end{align*}
By Theorem \ref{thm:presentacion minima}, it follows that $\cB(V)$ is presented by generators $x_1$, $x_2$, $x_3$, and relations
\begin{align*}
 x_1^2=x_2^2&=x_3^3= x_{\alpha_2+\alpha_3}^6=(\ad_cx_1)x_2=\left[x_{331}, x_{31}\right]_c = 0
 \\ [x_1,x_{332}]_c&=-q_{23}^{-1}\zeta^2[x_{132},x_3]_c+q_{13}x_3x_{132}.
\end{align*}

\noindent \vvii $\xymatrix{\ast^{-1}\ar@{-}[r]^{\zeta^2}&\bullet^{\zeta} \ar@{-}[r]^{-\zeta}&\circ^{-1}}$. We have the following positive roots:
\begin{align*}
\Delta_+^{\chi}=& \left\{\alpha_1,\alpha_2,\alpha_3,\alpha_2+\alpha_3,2\alpha_3+\alpha_1,\alpha_3+\alpha_1,
\alpha_2+\alpha_3+\alpha_1,\alpha_2+3\alpha_3+\alpha_1, \right.
\\& \left. \alpha_2+2\alpha_3+2\alpha_1,\alpha_2+2\alpha_3+\alpha_1,\alpha_2+3\alpha_3+2\alpha_1\right\}.
\end{align*}
Theorem \ref{thm:presentacion minima} implies that $\cB(V)$ has a presentation by generators $x_1$, $x_2$, $x_3$, and relations
$$ x_1^2=x_2^2=x_3^3= x_{\alpha_2+\alpha_3+\alpha_1}^6=(\ad_cx_3)^2x_2 =(\ad_cx_1)x_2= [x_{331},x_{31}]_c =0.$$

\noindent \vviii  $\xymatrix{ & \ast^{-1}\ar@{-}[rd]^{\zeta^2}&\\\circ^{-1} \ar@{-}[ru]^{-1}\ar@{-}[rr]^{-\zeta} & & \bullet^{\zeta} }$. For this braiding,
\begin{align*}
\Delta_+^{\chi}=& \left\{\alpha_1,\alpha_2,\alpha_3,\alpha_2+\alpha_1,\alpha_2+\alpha_3,\alpha_1+\alpha_3,
\alpha_1+2\alpha_3,\alpha_2+\alpha_1+\alpha_3, \right.
\\& \left. \alpha_2+\alpha_1+2\alpha_3,\alpha_2+2\alpha_1+2\alpha_3,\alpha_2+2\alpha_1+3\alpha_3\right\}.
\end{align*}
Following Theorem \ref{thm:presentacion minima}, $\cB(V)$ has a presentation by generators $x_1$, $x_2$, $x_3$, and relations
\begin{align*}
x_1^2&=x_2^2=x_3^3=x_{12}^2= x_{\alpha_2+\alpha_1+2\alpha_3}^6=(\ad_cx_3)^3x_2 = [x_{331},x_{31}]_c = 0
\\ x_{123}&=q_{23}(\zeta-\zeta^2)[x_{13},x_2]_c+q_{12}(1-\zeta)x_2x_{13}.
\end{align*}

\noindent \vix $\xymatrix{\ast^{-1}\ar@{-}[r]^{-1}&\circ^{-1}\ar@{-}[r]^{-\zeta^2}&\bullet^{\zeta^2}}$. Its positive roots are
\begin{align*}
\Delta_+^{\chi}=& \left\{\alpha_1,\alpha_2,\alpha_3,\alpha_2+\alpha_1,\alpha_1+\alpha_3,\alpha_1+2\alpha_3,
\alpha_2+\alpha_1+\alpha_3,\alpha_2+\alpha_1+2\alpha_3, \right.
\\& \left. \alpha_2+2\alpha_1+\alpha_3,\alpha_2+2\alpha_1+3\alpha_3,\alpha_2+2\alpha_1+2\alpha_3\right\}.
\end{align*}
According to Theorem \ref{thm:presentacion minima}, $\cB(V)$ is presented by generators $x_1$, $x_2$, $x_3$, and relations
\begin{align*}
x_1^2&=x_2^2=x_3^3= x_{\alpha_2+2\alpha_1+2\alpha_3}^6=(\ad_cx_2)x_3 =[x_{331},x_{31}]_c =[x_{3312},x_{312}]_c = 0
\\ x_{12}^2&= \left[ \left[x_{31},x_{312} \right]_c,x_1\right]_c=0.
\end{align*}
\end{example}

\medskip

\begin{example}{\textbf{Rank 3, row 18.}}\label{example:rank 3,row 18}
Let $\zeta\in\G'_9$ be primitive. We distinguish this Weyl
equivalence class of root systems because all the other unidentified
cases in rank $3$ have elements of $\G'_6$ labelling their vertices
and edges.

In this case, the roots in $\cO(\chi)$ are those such that $N_\alpha=9$, and $ \dim \cB(V)=9^93^4=3^{22}$.

\noindent \vi $\xymatrix{\circ^{\zeta}\ar@{-}[r]^{\zeta^8}&\circ^{\zeta} \ar@{-}[r]^{\zeta^8}&\circ^{\zeta^6}}$. We have that
\begin{align*}
\Delta_+^{\chi}=& \left\{\alpha_1,\alpha_2,\alpha_3,\alpha_1+\alpha_2,\alpha_2+\alpha_3,\alpha_1+\alpha_2+\alpha_3,
\alpha_1+\alpha_2+2\alpha_3,\alpha_1+2\alpha_2+3\alpha_3, \right.
\\& \left. \alpha_2+2\alpha_3,\alpha_1+2\alpha_2+4\alpha_3,\alpha_1+3\alpha_2+4\alpha_3,\alpha_1+2\alpha_2+2\alpha_3,
2\alpha_1+3\alpha_2+4\alpha_3\right\}.
\end{align*}
Following Theorem \ref{thm:presentacion minima}, $\cB(V)$ has a presentation by generators $x_1$, $x_2$, $x_3$, and relations
\begin{align*}
& x_3^3= (\ad_cx_1)^2x_2=(\ad_cx_1)x_3=(\ad_cx_2)^2x_1=(\ad_cx_2)^2x_3=0;
\\ & x_{\alpha}^9=0, \qquad \alpha\neq\alpha_3,\alpha_2+\alpha_3,\alpha_1+\alpha_2+\alpha_3,\alpha_1+2\alpha_2+3\alpha_3;
\\ & \left[ \left[x_{332},x_{3321} \right]_c,x_{32}\right]_c=0.
\end{align*}

\noindent \vii $\xymatrix{\circ^{\zeta}\ar@{-}[r]^{\zeta^8}&\circ^{\zeta^5}\ar@{-}[r]^{\zeta^4}&\circ^{\zeta^6}}$. For this braiding,
\begin{align*}
\Delta_+^{\chi}=& \left\{\alpha_1,\alpha_2,\alpha_3,\alpha_1+\alpha_2,\alpha_2+\alpha_3,\alpha_1+2\alpha_2,\alpha_2+2\alpha_3,
\alpha_1+\alpha_2+2\alpha_3, \right.
\\& \left. \alpha_1+\alpha_2+\alpha_3,\alpha_1+2\alpha_2+\alpha_3,\alpha_1+3\alpha_2+2\alpha_3,\alpha_1+2\alpha_2+2\alpha_3,
2\alpha_1+3\alpha_2+2\alpha_3 \right\}.
\end{align*}
Therefore $\cB(V)$ has a presentation by generators $x_1$, $x_2$, $x_3$, and relations
\begin{align*}
& x_3^3= (\ad_cx_1)^2x_2=(\ad_cx_1)x_3=(\ad_cx_2)^3x_1=(\ad_cx_2)^2x_3=0;
\\ & x_{\alpha}^9=0, \qquad \alpha\neq\alpha_3,\alpha_2+\alpha_3,\alpha_1+\alpha_2+\alpha_3,\alpha_1+2\alpha_2+\alpha_3;
\\ & (1+\zeta^4)\left[ \left[x_{123},x_{2} \right]_c,x_{3}\right]_c= q_{jk}\left[ \left[x_{123},x_{3} \right]_c,x_{2}\right]_c.
\end{align*}
\end{example}
\medskip

\begin{example}{\textbf{Rank 4, row 14.}}\label{example:rank 4,row 14}
Let $q\in\ku^{\times}$, $q\neq\pm1$. The roots in $\cO(\chi)$ are those such that $q_\alpha\in\{q,-q^{-1}\}$. Also, $\cB(V)$ is finite-dimensional iff $q$
has finite order. In such case, if $M=\ord q$, $N=\ord -q^{-1}$, we have $ \dim \cB(V)= 2^9M^3N^3$.

\noindent \vi $\xymatrix{\circ^{q}\ar@{-}[r]^{q^{-1}}&\circ^{q} \ar@{-}[r]^{q^{-1}}&\circ^{-1} \ar@{-}[r]^{-q}&\circ^{-q^{-1}}}$. For this braiding,
\begin{align*}
\Delta_+^{\chi}=& \left\{\alpha_1,\alpha_2,\alpha_3,\alpha_4,\alpha_1+\alpha_2,\alpha_2+\alpha_3,\alpha_3+\alpha_4,\alpha_1+\alpha_2+\alpha_3,
\alpha_2+\alpha_3+\alpha_4, \right.
\\& \left. \alpha_1+\alpha_2+\alpha_3+\alpha_4,\alpha_1+\alpha_2+2\alpha_3+\alpha_4,\alpha_1+2\alpha_2+3\alpha_3+\alpha_4,\right.
\\& \left.\alpha_2+2\alpha_3+\alpha_4,\alpha_1+2\alpha_2+2\alpha_3+\alpha_4,\alpha_1+2\alpha_2+3\alpha_3+2\alpha_4 \right\}.
\end{align*}
Following Theorem \ref{thm:presentacion minima}, $\cB(V)$ has a presentation by generators $x_1$, $x_2$, $x_3$, $x_4$, and relations
\begin{align*}
& (\ad_cx_1)^2x_2=(\ad_cx_2)^2x_1=(\ad_cx_2)^2x_3=(\ad_cx_4)^2x_3=0;
\\& x_3^2= (\ad_cx_1)x_3= (\ad_cx_1)x_4= (\ad_cx_2)x_4=0;
\\ & x_{\alpha}^M=0, \qquad \qquad \qquad \alpha= \alpha_1,\alpha_2,\alpha_1+\alpha_2+\alpha_3;
\\ & x_{\alpha}^N=0, \qquad \qquad \qquad \alpha= \alpha_4,\alpha_1+2\alpha_2+3\alpha_3+\alpha_4,\alpha_1+2\alpha_2+3\alpha_3+2\alpha_4.
\end{align*}

\noindent \vii $\xymatrix{\circ^{q}\ar@{-}[r]^{q^{-1}}&\circ^{-1}\ar@{-}[r]^{q}\ar@{-}[d]^{-1}&\circ^{-1}\ar@{-}[ld]^{-q^{-1}}\\&\circ^{-1}&}$. In this case,
\begin{align*}
\Delta_+^{\chi}=& \left\{\alpha_1,\alpha_2,\alpha_3,\alpha_4,\alpha_1+\alpha_2,\alpha_2+\alpha_3,\alpha_3+\alpha_4,\alpha_2+\alpha_4,
\alpha_2+\alpha_3+\alpha_4, \right.
\\& \left. \alpha_1+\alpha_2+\alpha_3,\alpha_1+\alpha_2+\alpha_4,\alpha_1+2\alpha_2+\alpha_3+\alpha_4,\right.
\\& \left.\alpha_1+2\alpha_2+\alpha_4,\alpha_1+2\alpha_2+\alpha_3+2\alpha_4,\alpha_1+\alpha_2+\alpha_3+\alpha_4 \right\}.
\end{align*}
Therefore $\cB(V)$ has a presentation by generators $x_1$, $x_2$, $x_3$, $x_4$, and relations
\begin{align*}
& (\ad_cx_1)^2x_2=(\ad_cx_1)x_3=(\ad_cx_1)x_4= \left[x_{123},x_2\right]_c= x_2^2=x_3^2=x_4^2=x_{24}^2=0;
\\ & x_{234}+\frac{1}{2}(1+q)q_{43}\left[x_{24},x_3\right]_c-q_{23}(1+q^{-1})x_3x_{24}=0;
\\ & x_{\alpha}^M=0, \qquad \alpha= \alpha_1,\alpha_2+\alpha_3,\alpha_1+\alpha_2+\alpha_3;
\\ & x_{\alpha}^N=0, \qquad \alpha= \alpha_3+\alpha_4,\alpha_1+2\alpha_2+\alpha_4,\alpha_1+2\alpha_2+\alpha_3+2\alpha_4.
\end{align*}

\noindent \viii $\xymatrix{\circ^{q}\ar@{-}[r]^{q^{-1}}&\circ^{-1} \ar@{-}[r]^{-1}&\circ^{-1} \ar@{-}[r]^{-q}&\circ^{-q^{-1}}}$. For this braiding,
\begin{align*}
\Delta_+^{\chi}=& \left\{\alpha_1,\alpha_2,\alpha_3,\alpha_4,\alpha_1+\alpha_2,\alpha_2+\alpha_3,\alpha_3+\alpha_4,\alpha_1+\alpha_2+\alpha_3,
\alpha_2+\alpha_3+\alpha_4, \right.
\\& \left. \alpha_1+2\alpha_2+\alpha_3,\alpha_1+\alpha_2+\alpha_3+\alpha_4,\alpha_1+2\alpha_2+2\alpha_3+\alpha_4,\right.
\\& \left.\alpha_2+2\alpha_3+\alpha_4,\alpha_1+2\alpha_2+\alpha_3+\alpha_4,\alpha_1+\alpha_2+2\alpha_3+\alpha_4 \right\}.
\end{align*}
Then $\cB(V)$ has a presentation by generators $x_1$, $x_2$, $x_3$, $x_4$, and relations
\begin{align*}
& x_2^2=x_3^2=x_{23}^2= (\ad_cx_i)^{1+m_{ij}}x_j=0, \qquad \qquad i\neq 2,3;
\\ & \big[ \left[ x_{1234} , x_2 \right]_c , x_3 \big]_c = q_{23}(1+q^{-1})\big[ \left[ x_{1234} , x_3 \right]_c , x_2 \big]_c;
\\ & x_{\alpha}^M=0, \qquad \qquad \qquad \alpha= \alpha_1,\alpha_2+2\alpha_3+\alpha_4,\alpha_1+\alpha_2+2\alpha_3+\alpha_4;
\\ & x_{\alpha}^N=0, \qquad \qquad \qquad \alpha= \alpha_4,\alpha_1+2\alpha_2+\alpha_3,\alpha_1+2\alpha_2+\alpha_3+\alpha_4.
\end{align*}
\end{example}

\begin{example}{\textbf{Rank 4, row 22.}}\label{example:rank 4,row 22}
Let $\zeta\in\G'_4$. This Weyl equivalence class contains eight
different diagrams:
\begin{align*}
&\mathbf{(i)}
\xymatrix@=.6cm{\circ^{-\zeta}\ar@{-}[r]^{\zeta}&\circ^{-1}
\ar@{-}[r]^{-\zeta}&\circ^{\zeta} \ar@{-}[r]^{\zeta}&\circ^{-\zeta}}
& \mathbf{(ii)}\xymatrix@=.6cm{\circ^{-1}\ar@{-}[r]^{-\zeta}&\circ^{-1}
\ar@{-}[r]^{\zeta}&\circ^{-1} \ar@{-}[r]^{\zeta}&\circ^{-\zeta}}
\\ &\mathbf{(iii)} \xymatrix@=.6cm{\circ^{-\zeta}\ar@{-}[r]^{\zeta}&\circ^{-1}\ar@{-}[r]^{-\zeta}
\ar@{-}[d]_{-1}&\circ^{\zeta}\ar@{-}[ld]^{-\zeta}\\&\circ^{-1}&}
& \mathbf{(iv)}\xymatrix@=.6cm{\circ^{-1}\ar@{-}[r]^{-\zeta}&\circ^{\zeta}\ar@{-}[r]^{-\zeta}
\ar@{-}[d]_{-1}&\circ^{-1}\ar@{-}[ld]^{-\zeta}\\&\circ^{-1}&}
\\ &\mathbf{(v)} \xymatrix@=.6cm{\circ^{-1}\ar@{-}[r]^{\zeta}&\circ^{-\zeta}
\ar@{-}[r]^{\zeta}&\circ^{-1} \ar@{-}[r]^{\zeta}&\circ^{-\zeta}}
& \mathbf{(vi)}\xymatrix@=.6cm{\circ^{-1}\ar@{-}[r]^{\zeta}&\circ^{-1}\ar@{-}[r]^{-\zeta}
\ar@{-}[d]_{-1}&\circ^{-1}\ar@{-}[ld]^{-\zeta}\\&\circ^{-1}&}
\\ &\mathbf{(vii)}
\xymatrix@=.6cm{\circ^{-1}\ar@{-}[r]^{-\zeta}&\circ^{\zeta}
\ar@{-}[r]^{-1}&\circ^{-1} \ar@{-}[r]^{\zeta}&\circ^{-\zeta}}
& \mathbf{(viii)}\xymatrix@=.6cm{\circ^{-1}\ar@{-}[r]^{\zeta}&\circ^{-1}
\ar@{-}[r]^{-1}&\circ^{-1} \ar@{-}[r]^{\zeta}&\circ^{-\zeta}}
\end{align*}
Each ones of the associated Nichols algebras have dimension $2^{42}$. They are presented by
generators $x_1$, $x_2$, $x_3$, $x_4$, and relations:
\begin{itemize}
  \item if $N_{\alpha}=4$:  $\quad x_{\alpha}^{N_{\alpha}}$;
  \item if $q_{ii}=-1$: $\quad x_i^2$;
  \item if $q_{ii}=q_{jj}=\widetilde{q_{ij}}=-1$ (diagrams \textbf{(iii)}, \textbf{(vi)}, \textbf{(viii)}): $\quad x_{ij}^2$;
  \item if $\widetilde{q_{ij}}=1$: $\quad (\ad_cx_i)x_j$;
  \item if $\widetilde{q_{ij}}=q_{ii}^{-1}=\pm\zeta$: $\quad (\ad_cx_i)^2x_j$;
  \item if $\widetilde{q_{ij}}=-1$, $q_{ii}=\zeta$ (diagrams \textbf{(iv)}, \textbf{(vii)}): $\quad \left[x_{iij},x_{ij} \right]_c$;
  \item if $\widetilde{q_{ij}}=\widetilde{q_{kj}}=-\zeta$, $\widetilde{q_{ik}}=-1$
  (diagrams \textbf{(iii)}, \textbf{(iv)}, \textbf{(vi)}):
  $$ x_{ijk}-\zeta q_{jk}[x_{ik},x_j]_c-q_{ij}(1+\zeta)x_jx_{ik} ; $$
  \item if $-q_{jj}=\widetilde{q_{ij}}\widetilde{q_{kj}}=\widetilde{q_{ik}}=1$ (diagrams \textbf{(i)}, \textbf{(ii)}, \textbf{(iii)}, \textbf{(vi)}): $\quad [x_{ijk},x_j]_c$;
  \item if $q_{jj}=\widetilde{q_{ij}}=-\widetilde{q_{kj}}=\zeta$, $\widetilde{q_{ik}}=1$ (diagram \textbf{(i)}):
  $\quad \left[\left[\left[x_{ijk},x_j\right]_c,x_j\right]_c,x_j\right]_c$;
  \item if $q_{ii}=q_{jj}-1$, $-\widetilde{q_{ij}}=\widetilde{q_{kj}}=\zeta$, $\widetilde{q_{ik}}=1$
  (diagram \textbf{(ii)}): $\quad \left[\left[x_{ij}, \left[x_{ij},x_{ijk} \right]_c\right]_c, x_j\right]_c$;
  \item for diagram \textbf{(v)}: $\left[\left[ x_{123}, \left[ x_{1234} , x_3 \right]_c\right]_c, x_{23} \right]_c$.
\end{itemize}
\end{example}

\medskip

\section{The other unidentified Nichols algebras}

In this Section we will consider the remaining braidings of
unidentified type. We will consider four differents subfamilies,
according with the shape of their generalized Dynkin diagrams. The
first subfamily that we will study is closely related with diagrams
of type $D_5$, $E_6$, $E_7$, but they contain non-Cartan vertices
labeled with $-1$ and small orders on the $q_{ii}$ (in fact, they
are roots of unity of order 3 or 4).
\medskip

\begin{theorem}\label{thm:casos similar a E}
Let $(V,c)$ a braided vector space of diagonal type of dimension $\theta\in\{5,6,7\}$, whose generalized Dynkin
diagram belongs to rows 11, 14, 17, 19 or 21 of \cite[Table 4]{H-classif RS}. Then $\cB(V)$ is presented by
generators $x_1,\ldots,x_{\theta}$ and relations
\begin{itemize}
  \item[$\odot$] if $N_{\alpha}\neq 2$ or $N_\alpha=2$, $\alpha=\alpha_i$: $\quad  x_{\alpha}^{N_{\alpha}}$;
  \item[$\odot$] if $i,j$ are such that $q_{ii}=\widetilde{q_{ij}}=q_{jj}=-1$: $\quad x_{ij}^2$;
  \item[$\odot$] if $i,j$ are such that $\widetilde{q_{ij}}=1$ (respectively, $\widetilde{q_{ij}}=q_{ii}^{-1}\neq-1$):
  $$ (\ad_cx_i)x_j \qquad (\mbox{resp. }(\ad_cx_i)^2x_j; $$
  \item[$\odot$] if $i,j,k$ are such that $\widetilde{q_{ik}}=1$, $q_{jj}=-1$, $\widetilde{q_{ij}}\widetilde{q_{jk}}=1$:
  $\quad [x_{ijk},x_j]_c;$
  \item[$\odot$] if $i,j,k$ are such that $\widetilde{q_{ik}}=\widetilde{q_{ij}}=\widetilde{q_{jk}}=:\zeta\in\G'_3$:
  $$ x_{ijk}-q_{jk}\zeta^2[x_{ik},x_j]_c-q_{ij}(1-\zeta)x_jx_{ik}; $$
  \item[$\odot$] if $i,j,k$ are such that $\widetilde{q_{ik}}=\widetilde{q_{jk}}=:\zeta\in\G'_4$, $\widetilde{q_{ij}}=-1$:
  $$ x_{ijk}+q_{jk}\zeta[x_{ik},x_j]_c-q_{ij}(1-\zeta)x_jx_{ik}. $$
\end{itemize}
\end{theorem}
\pf
First, note that $a_{ij}\in\{0,-1\}$ for all these braidings; that is, for each pair $i\neq j$, either $\widetilde{q_{ij}}=1$,
or else $\widetilde{q_{ij}}\neq-1$, $q_{ii}\in\{-1,\widetilde{q_{ij}}^{-1}\}$. This explains why we only consider quantum Serre
relations associated to $a_{ij}=0,1$. We derive also that the Cartan vertices are those labeled with $q_{ii}\neq-1$, which
explains why we need only the power root vectors associated to simple roots, or to other roots such that $N_\alpha\neq2$.

Finally we look at the other needed relations. As the $a_{ij}$'s take only two values, we need few extra relations.
\epf

\medskip

The second subfamily seems close to type super $D(n)$, but with certain 'degeneration' and small order on the labels
of the vertices (roots of unity of order 2, 3, 5). With these diagrams and the corresponding to the previous family
we cover all the cases in rank $5,6,7$.
\medskip

\begin{theorem}\label{thm:casos similar a C}
Let $(V,c)$ a braided vector space of diagonal type of dimension $\theta\in\{3,4,5,6\}$, whose generalized Dynkin
diagram belongs to row 15 of \cite[Table 2]{H-classif RS}, or row 18 of \cite[Table 3]{H-classif RS}, or rows 12,
13, 15 or 18 of \cite[Table 4]{H-classif RS}. Then $\cB(V)$ is presented by
generators $x_1,\ldots,x_{\theta}$ and relations
\begin{itemize}
  \item[$\odot$] if $N_{\alpha}\neq 2$: $\quad  x_{\alpha}^{N_{\alpha}}$;
  \item[$\odot$] if $q_{ii}=-1$: $\qquad x_i^2$;
  \item[$\odot$] if $i,j$ are such that $\widetilde{q_{ij}}=1$: $\qquad (\ad_cx_i)x_j;$
  \item[$\odot$] if $i,j$ are such that $\widetilde{q_{ij}}=q_{ii}^{-1}\neq-1$: $\qquad (\ad_cx_i)^2x_j;$
  \item[$\odot$] if $i,j$ are such that $\widetilde{q_{ij}}=q_{ii}^{-2}\in\G'_5$: $\qquad (\ad_cx_i)^3x_j;$
  \item[$\odot$] if $i,j,k$ are such that $\widetilde{q_{ik}}=1$, $q_{jj}=-1$, $\widetilde{q_{ij}}\widetilde{q_{jk}}=1$:
  $\quad [x_{ijk},x_j]_c;$
  \item[$\odot$] if $i,j,k$ are such that $\widetilde{q_{ik}}=\widetilde{q_{ij}}=\widetilde{q_{jk}}=:\zeta\in\G'_3$:
  $$ x_{ijk}-q_{jk}\zeta^2[x_{ik},x_j]_c-q_{ij}(1-\zeta)x_jx_{ik}; $$
  \item[$\odot$] if $i,j,k$ are such that $\widetilde{q_{ij}}=:\zeta\in\G'_5$, $\widetilde{q_{ik}}=\widetilde{q_{jk}}=\zeta^2$:
  $$ x_{ijk}-q_{jk}\zeta^3[x_{ik},x_j]_c-q_{ij}(1-\zeta^2)x_jx_{ik}. $$
  \item[$\odot$] if $i,j,k$ are such that $\widetilde{q_{ik}}=1$, $q_{jj}=\widetilde{q_{ij}}=\widetilde{q_{jk}}^2\in\G'_3$:
  $$ \left[ \left[ x_{ijk}, x_j\right]_c, x_j\right]_c; $$
  \item[$\odot$] if $i,j,k$ are such that $\widetilde{q_{ik}}=1$, $q_{ii}=q_{jj}=-1$, $\widetilde{q_{jk}}=\widetilde{q_{ik}}^{-2}$:
  $$ \left[ \left[ x_{ij},x_{ijk}\right]_c, x_j\right]_c. $$
\end{itemize}
\end{theorem}
\pf
For these braidings, $a_{ij}\in\{0,-1,-2\}$. Moreover, the Cartan vertices are those such that $q_{ii}\neq-1$. In fact, when
$q_{ii}\in\G'_3$ and there exists $j$ such that $a_{ij}=-2$, then $q_{ij}=q_{ii}$, so $i$ is a Cartan vertex (but we do not need
the corresponding quantum Serre relation). Therefore we need only the power root vectors associated to a simple root, or to
other roots such that $N_\alpha\neq2$.

The remaining relations we need to generate the ideal expresses the similarity with the super $D(n)$ case.
\epf

\medskip

The following subfamily keeps certain similarity with diagrams of
type super $F(4)$ for small orders on the labels of the vertices (3,
6).
\medskip

\begin{theorem}\label{thm:casos similar a F(4)}
Let $(V,c)$ a braided vector space of diagonal type of dimension $\theta=4$, whose generalized Dynkin
diagram belongs to rows 20 or 21 of \cite[Table 3]{H-classif RS}. Then $\cB(V)$ is presented by
generators $x_1,x_2,x_3,x_4$ and relations
\begin{itemize}
  \item[$\odot$] if $N_{\alpha}=3,6$: $\quad  x_{\alpha}^{N_{\alpha}}$;
  \item[$\odot$] if $q_{ii}=-1$: $\qquad x_i^2$;
  \item[$\odot$] if $i,j$ are such that $\widetilde{q_{ij}}=1$: $\qquad (\ad_cx_i)x_j;$
  \item[$\odot$] if $i,j$ are such that $\widetilde{q_{ij}}=q_{ii}^{-1}\neq-1$: $\qquad (\ad_cx_i)^2x_j;$
  \item[$\odot$] if $i,j$ are such that $\widetilde{q_{ij}}=q_{ii}^{-2}$, $q_{ii}\in\G'_6$: $\qquad (\ad_cx_i)^3x_j;$
  \item[$\odot$] if $i,j,k$ are such that $q_{jj}=-1$, $\widetilde{q_{ij}}\widetilde{q_{jk}}=1$, $\widetilde{q_{ik}}=1$:
  $\quad [x_{ijk},x_j]_c;$
  \item[$\odot$] if $i,j,k$ are such that $\widetilde{q_{ik}}=1$, $q_{jj}=\widetilde{q_{ij}}=\widetilde{q_{jk}}^2\in\G'_3$:
  $$ \left[ \left[ x_{ijk}, x_j\right]_c, x_j\right]_c; $$
  \item[$\odot$] if $i,j,k$ are such that $\widetilde{q_{ik}}=1$, $q_{ii}=q_{jj}=-1$, $\widetilde{q_{jk}}=\widetilde{q_{ik}}\in\G'_3$:
  $$ \left[ \left[ x_{ij},x_{ijk}\right]_c, x_j\right]_c; $$
  \item[$\odot$] if $i,j,k$ are such that $\widetilde{q_{ik}}=\widetilde{q_{ij}}=\widetilde{q_{jk}}=:\zeta\in\G'_3$:
  $$ x_{ijk}-q_{jk}\zeta^2[x_{ik},x_j]_c-q_{ij}(1-\zeta)x_jx_{ik}; $$
  \item[$\odot$] if $i,j,k$ are such that $\widetilde{q_{ik}}=1$, $-q_{jj}=\widetilde{q_{ij}}=\widetilde{q_{jk}}=:\zeta\in\G'_3$:
  $$ [x_i,x_{jjk}]_c+q_{jk}[x_{ijk},x_j]_c+q_{ij}(\zeta-\zeta^2)x_jx_{ijk}. $$
\end{itemize}
\end{theorem}
\pf
It is analogous to the proof of Theorem \ref{thm:casos similar a C}.
\epf

\begin{remark}
If the diagram of $(V,c)$ belongs to row 20 of \cite[Table 3]{H-classif RS}, then
$$\dim \cB(V)=2^{13}3^{10}.$$
On the other hand, if it belongs to row 21 of \cite[Table 3]{H-classif RS}, then
$$\dim \cB(V)=2^{20}3^{16}.$$
\end{remark}

\medskip

The last subfamily contains only two Weyl equivalence classes, where
the labels of the vertices are roots of unity of order 2, 3, 6.
\medskip

\begin{theorem}\label{thm:casos restantes}
Let $(V,c)$ a braided vector space of diagonal type of dimension $\theta=3$ (resp. $\theta=4$), whose generalized Dynkin
diagram belongs to row 16 of \cite[Table 2]{H-classif RS},  (resp. row 17 of \cite[Table 3]{H-classif RS}). Then $\cB(V)$ is presented by
generators $x_1,\ldots,x_{\theta}$ and relations
\begin{itemize}
  \item[$\odot$] if $N_{\alpha}=6$: $\quad  x_{\alpha}^6$;
  \item[$\odot$] if $N_i=2,3$: $\qquad x_i^{N_i}$;
  \item[$\odot$] if $i,j$ are such that $\widetilde{q_{ij}}=1$: $\qquad (\ad_cx_i)x_j;$
  \item[$\odot$] if $i,j$ are such that $\widetilde{q_{ij}}=q_{ii}^{-1}\neq-1$: $\qquad (\ad_cx_i)^2x_j;$
  \item[$\odot$] if $i,j$ are such that $\widetilde{q_{ij}}=q_{ii}=q_{jj}=-1$: $\qquad x_{ij}^2;$
  \item[$\odot$] if $i,j$ are such that $q_{ii}\in\G'_3$, $q_{jj}=\widetilde{q_{ij}}=-1$:
  $\quad [x_{iij},x_{ij}]_c;$
  \item[$\odot$] if $i,j,k$ are such that $\widetilde{q_{ij}}=-\widetilde{q_{ik}}=-\widetilde{q_{jk}}=:\zeta\in\G'_3$:
  $$ x_{ijk}+q_{jk}\zeta^2[x_{ik},x_j]_c+q_{ij}\zeta^2x_jx_{ik}; $$
  \item[$\odot$] if $i,j,k$ are such that $\widetilde{q_{ij}}=-1$, $\widetilde{q_{ik}}^{-1}=-\widetilde{q_{jk}}=:\zeta\in\G'_3$:
  $$ x_{ijk}+\frac{1}{3}q_{jk}\zeta[x_{ik},x_j]_c+q_{ij}\zeta^2x_jx_{ik}; $$
  \item[$\odot$] if $i,j,k,l$ are such that $\widetilde{q_{ij}}\in\G_6'$, $q_{jj}^{-1}=\widetilde{q_{kj}}=\widetilde{q_{ij}}^2$, $\widetilde{q_{kl}}=q_{kk}=-1$, $\widetilde{q_{il}}=\widetilde{q_{ik}}=\widetilde{q_{jl}}=1$:
  $$\left[\left[ x_{ijk}, \left[ x_{ijkl} , x_k \right]_c\right]_c, x_{jk} \right]_c.$$
\end{itemize}
\end{theorem}
\pf Now $a_{ij}\in\{0,-1,-2\}$, and the Cartan vertices are those
such that $q_{ii}^2, q_{ii}^3\neq1$. The proof follows then as in
the previous cases. \epf

\begin{remark}
If the diagram of $(V,c)$ belongs to \cite[Table 2, row 16]{H-classif RS}, $\dim \cB(V)=2^73^6$.

If it belongs to \cite[Table 3, row 17]{H-classif RS}, then $\dim \cB(V)=2^{19}3^{15}$.
\end{remark}

\subsection*{Acknowledgements}

The author thanks Leandro Vendram\'in for all the explanations about
SARNA program \cite{GHV} to obtain root systems. He thanks also
Hiroyuki Yamane for many helpful discussions about the presentation.


\begin{thebibliography}{XXXX}

\bibitem[AA]{AA} N. Andruskiewitsch and I. Angiono, \emph{On Nichols algebras with generic braiding}, Modules and Comodules.
Trends in Mathematics. Brzezinski, T.; G\'omez Pardo, J.L.; Shestakov, I.; Smith, P.F. (Eds.), pp. 47--64 (2008).

\bibitem[AS1]{AS lifting} N. Andruskiewitsch and H.-J. Schneider,
\emph{Lifting of quantum linear spaces and pointed Hopf algebras of order
$p^3$}, J. Algebra \textbf{209}, 658--691 (1998).

\bibitem[AS2]{AS Fin-qg} \bysame, \emph{Finite quantum groups and Cartan matrices}, Adv. Math. \textbf{154}, 1--45 (2000).

\bibitem[AS3]{AS Pointed HA} \bysame, \emph{Pointed Hopf algebras}, ``New directions in Hopf algebras'',
MSRI series Cambridge Univ. Press; 1--68 (2002).

\bibitem[AS4]{AS Class} \bysame, \emph{On the classification of finite-dimensional pointed Hopf algebras},
Annals of Mathematics Vol. \textbf{171} (2010), No. 1, 375--417.

\bibitem[A1]{A-Standard} I. Angiono, \emph{On Nichols algebras with standard braiding}, Algebra and Number Theory Vol. 3,
No. \textbf{1}, 35--106, (2009).

\bibitem[A2]{A-presentation} I. Angiono, \emph{A presentation by generators and relations of Nichols algebras of diagonal type and convex orders on root
systems}, submitted. arXiv:1008.4144.

\bibitem[A3]{A-exp presentation} \bysame, \emph{On Nichols algebras of diagonal type}, submitted. arXiv:1104.0268.

\bibitem[AAY]{AAY} N. Andruskiewitsch, I. Angiono and H. Yamane, \emph{On pointed Hopf superalgebras}.
Contemp. Math. \textbf{544}, to appear. arxiv:1009.5148, 18 pp.

\bibitem[CH]{CH1} M. Cuntz and I. Heckenberger, \emph{Weyl groupoids with at most three objects}. J. Pure Appl. Algebra \textbf{213}, No. 6 (2009),
1112--1128.

\bibitem[GHV]{GHV} M. Gra\~na, I. Heckenberger and L. Vendram\'in, \emph{SARNA: Arithmetic Root Systems and Nichols Algebras
Software}, a program to calculate the Arithmetic Root System associated to a braid of diagonal type. Available at \texttt{http://mate.dm.uba.ar/~lvendram/}

\bibitem[H1]{H-Weyl gpd} I. Heckenberger, \emph{The Weyl groupoid of a Nichols algebra of diagonal type}, Invent. Math. \textbf{164}, 175--188 (2006).

\bibitem[H2]{H-classif RS} \bysame, \emph{Classification of arithmetic root systems}, Adv. Math. \textbf{220} (2009) 59--124.

\bibitem[HS]{HS} I. Heckenberger and H.-J. Schneider, \emph{Right coideal subalgebras of Nichols algebras and the Duflo order on the Weyl groupoid}.
\texttt{math.QA/0909.0293}.

\bibitem[HY]{HY} I. Heckenberger and H. Yamane, \emph{A generalization of Coxeter groups, root systems, and Matsumoto's theorem}, Math. Z.
\textbf{259} (2008), 255--276.

\bibitem[Kh]{Kh} V. Kharchenko, \emph{A quantum analog of the Poincare-Birkhoff-Witt theorem}, Algebra and Logic \textbf{38}, (1999), 259-276.

\bibitem[L]{L} G. Lusztig, \emph{Introduction to quantum groups}, Birkh\"auser (1993).


\bibitem[U]{U} S. Ufer, \emph{PBW bases for a class of braided Hopf algebras}, J. Alg. \textbf{280} (2004) 84-119.

\bibitem[Y]{Y} H. Yamane, \emph{Quantized enveloping algebras associated to simple Lie superalgebras and their universal $R$-matrices}
Publ. Res. Inst. Math. Sci.  \textbf{30} (1994), 15--87.
\end{thebibliography}
\end{document}